\documentclass[12pt, a4paper]{article}

\usepackage[utf8]{inputenc}
\usepackage{titlesec}
\usepackage[english]{babel}
\usepackage{graphicx}
\usepackage{float}
\usepackage{fancyhdr}
\usepackage[matrix,arrow,curve]{xy}
\usepackage{pstricks} 
\usepackage{amsmath,amsfonts,verbatim,afterpage,theorem,euscript,mathrsfs,amssymb}
\usepackage{amsfonts}
\usepackage{amssymb}
\usepackage{array}
\usepackage{dsfont}
\usepackage{hyperref}
\usepackage{authblk}

\newtheorem{Definition}{Definition}[section]
\newtheorem{Proposition}{Proposition}[section]

\newtheorem{Lemma}{Lemma}[section]
\newtheorem{Theorem}{Theorem}
\newtheorem{Corollary}{Corollary}[section]
\newtheorem{Remark}{Remark}[section]

\def \vu{\textbf{u}}
\def \vb{\textbf{b}}
\def \vv{\textbf{v}}
\def \vc{\textbf{c}}
\def \vd{\textbf{d}}
\def \vw{\textbf{w}}

\def \Ri{\mathcal{R}}

\newcommand \vue{\textbf{u}_\varepsilon}
\newcommand \vbe{\textbf{b}_\varepsilon}

\newcommand \pe{p_\varepsilon}
\newcommand \qe{q_\varepsilon}

\newcommand \te{\theta_\varepsilon}

\newcommand\vN{{\bf  \nabla}}

\def \P{\mathbb{P}}
\def \F{\mathbb{F}}
\def \G{\mathbb{G}}

\def \Rt{\mathbb{R}^{3}}

\def \pv{{\bf{Proof.}}~}

\def \ds{\displaystyle}
\newcommand \Endproof{\hfill $\diamond$}

\begin{document}

\title{Discretely self-similar solutions for 3D MHD equations and global weak solutions in weighted $L^2$ spaces}
\author{Pedro Gabriel Fern\'andez-Dalgo\footnote{LaMME, Univ Evry, CNRS, Universit\'e Paris-Saclay, 91025, Evry, France } \footnote{e-mail : pedro.fernandez@univ-evry.fr} , Oscar Jarr\'in \footnote{Direcci\'on de Investigaci\'on y Desarrollo (DIDE), Universidad T\'ecnica de Ambato, Ambato, Ecuador} \footnote{e-mail : or.jarrin@uta.edu.ec}}
\date{}\maketitle

\maketitle

\begin{abstract}
This paper deals with the existence of global weak solutions for 3D MHD equations when the initial data belong to the weighted spaces  $L^2_{w_\gamma}$, with $w_\gamma(x)=(1+\vert x\vert)^{-\gamma}$ and $0 \leq \gamma \leq 2$. Moreover, we prove the existence of discretely self-similar solutions for
3D MHD equations for discretely self-similar initial data which are locally square integrable. Our methods are inspired  of a recent work  \cite{PF_PG} for the Navier-Stokes equations. 
\end{abstract}
 
\noindent{\bf Keywords : } MHD equations, weighted $L^2$ spaces, discretely self-similar solutions, energy controls.\\
\noindent{\bf AMS classification : }  35Q30, 76D05.

\section{Introduction}

The Cauchy problem for the incompressible and homogeneous magneto-hydrodynamic equations (MHD) equations in the whole space $\Rt$ writes down as: 
\begin{equation}\label{MHD}
(\text{MHD}) \left\{ \begin{array}{ll}\vspace{2mm} 
\partial_t \vu = \Delta \vu - (\vu \cdot \nabla) \vu + (\vb \cdot \nabla) \vb - \nabla p + \vN \cdot \F,  \\ \vspace{2mm}
\partial_{t} \vb = \Delta \vb - (\vu \cdot \nabla) \vb + (\vb \cdot \nabla) \vu 
, \\ \vspace{2mm}
\vN \cdot  \vu=0, \,  \vN \cdot \vb =0, \\ \vspace{2mm}
\vu(0,\cdot)=\vu_0, \, \vb(0,\cdot)=\vb_0,  
\end{array}
\right.     
\end{equation}
where the fluid velocity field $\vu:[0,+\infty)\times \Rt \to \Rt$, the magnetic field $\vb:[0,+\infty)\times \Rt \to \Rt$  and the fluid pressure $p:[0,+\infty)\times \Rt \to \mathbb{R}$ are the unknowns, and the fluid velocity at $t=0$:  $\vu_0:\Rt\to \Rt$, the magnetic field at $t=0$: $\vb_0:\Rt\to \Rt$, and the tensor $\F=(F_{i,j})_{1\leq i,j\leq 3}$ (where $F_{i,j}: [0,+\infty)\times \Rt \to \mathbb{R}^9$) whose divergence $\vN \cdot F$ represents a volume force applied to the fluid, are the data of the problem. \\

In this article, we will focus on the following simple generalisation of (MHD) equations:
\begin{equation}\label{MHDG}
(\text{MHDG}) \left\{ \begin{array}{ll}\vspace{2mm} 
\partial_t \vu = \Delta \vu - (\vu \cdot \nabla) \vu + (\vb \cdot \nabla) \vb - \nabla p + \vN \cdot \F,  \\ \vspace{2mm}
\partial_{t} \vb = \Delta \vb - (\vu \cdot \nabla) \vb + (\vb \cdot \nabla) \vu -\vN q + 
\vN  \cdot \G 
, \\ \vspace{2mm}
\vN \cdot  \vu=0, \,  \vN \cdot \vb =0, \\ \vspace{2mm}
\vu(0,\cdot)=\vu_0, \, \vb(0,\cdot)=\vb_0,  
\end{array}
\right.     
\end{equation}
where in the second equation we have added an extra gradient term $\vN q$, which is an unknown, and an extra tensor field $\G=(G_{i,j})_{1\leq i,j \leq 3}$ which is a datum. This generalized system   does not present extra mathematical difficulties but it appears in physical models when Maxwell's displacement currents are considered \cite{Bat,Shercliff}. Moreover, we construct solutions for (MHDG) such that $G=0$ implies $q=0$ (see the equation (\ref{q}) below), and it justifies the fact that (MHDG) generalizes (MHD) from the mathematical point of view.\\

In the recent work \cite{PF_PG} due to P. Fernandez \& P.G. Lemarié-Rieusset, which deals with the homogeneous and incompressible Navier-Stokes equations  in the whole space $\Rt$:
\begin{equation*}
(\text{NS}) \left\{ \begin{array}{ll}\vspace{2mm} 
\partial_t \vu = \Delta \vu - (\vu \cdot \nabla) \vu  - \nabla p + \vN \cdot \F,  \\ \vspace{2mm}
\vN \cdot  \vu=0, \quad \vu(0,\cdot)=\vu_0, 
\end{array}
\right.     
\end{equation*}
the authors established \emph{new energy controls} which have two prominent applications: the first one is to develop a theory to construct infinite-energy global weak solutions for (NS) arising from  \emph{initial datum} $\vu_0$  belonging to the weighted space $L^{2}((1+\vert x \vert)^{-\gamma}dx)$, where $0 < \gamma \leq 2$, and the second one is to give a new proof of the existence of discretely self-similar solutions for discretely self-similar initial data which are locally square integrable (proven before in \cite{CW18} by Chae and Wolf and in \cite{BT19} by Bradshaw and Tsai). \\ 

In \cite{BTK}, Bradshaw, Tsai \& Kukavika give an improvement of the existence theorem in \cite{PF_PG} with respect to the space of initial data. For other constructions of  infinite-energy weak solutions for the (NS) equations see the articles \cite{Ba06,BTs19,JS14,KS07,LR99} and the books \cite{LR02,LR16}. \\

 The main purpose of this article is to adapt the energy methods given in \cite{PF_PG} for (NS) to the more general setting of the coupled  system (MHDG). We remark (Proposition \ref{Prop-approx-leray})  that approximate solutions for (NS) and (MHDG) admit an energy balance which have a similar structure.
 Our first result reads as follows: 
\begin{Theorem} \label{weightedMHD} Let $ 0 \leq \gamma\leq 2$. Let  $\vu_{0}, \vb_{0}$ be  divergence-free vector fields such that $\vu_{0}, \vb_{0}\in L^2_{w_\gamma}(\mathbb{R}^3)$. Let    $\mathbb{F}$ and $\G$ be    tensors  such that $\mathbb{F}, \G  \in L^2((0,+ \infty), L^2_{w_\gamma})$.  Then, the system (MHDG)
has a solution $(\vu, \vb,p,q)$  which satisfies :
 \begin{itemize} 
 \item[$\bullet$]  $\vu, \vb$ belong to $L^\infty((0,T), L^2_{w_\gamma})$ and $\vN\vu$, $\vN \vb$ belong to $L^2((0,T),L^2_{w_\gamma})$, for all $0 < T < + \infty$.
 \item[$\bullet$] The pressure $p$ and the term $q$ are related to $\vu$, $\vb$, $\mathbb{F}$ and $\G$ by 
 $$ p =\sum_{1 \leq i,j\leq 3} \Ri_i \Ri_j(u_iu_j-b_i b_j -F_{i,j})$$
 and
 \begin{equation}\label{q} 
  q =-\sum_{1 \leq i,j\leq 3} \Ri_i \Ri_j(G_{i,j}).    
 \end{equation}
 \item[$\bullet$] The map $t\in [0,+\infty)\mapsto (\vu(t), \vb(t ))$ is weakly continuous from $[0,+\infty)$ to $L^{2}_{w_\gamma}$, and is strongly continuous at $t=0$ :
 $$ \lim_{t\rightarrow 0}    \Vert  ( \vu(t,\cdot) - \vu_{0} , \vb(t,\cdot) - \vb_{0} ) \Vert_{L^{2}_{w_\gamma}} = 0.$$
 \item[$\bullet$] the solution $(\vu, \vb, p,q)$ is suitable : there exist a non-negative locally finite measure $\mu$ on $(0,+\infty)\times\mathbb{R}^3$ such that
\begin{equation}\label{energloc1}
\begin{split}
   \partial_t(\frac {\vert\vu \vert^2 
    + |\vb|^2 }2)=&\Delta(\frac {\vert\vu \vert^2 + |\vb|^2 }2)-\vert\vN \vu \vert^2 - |\vN \vb|^2 \\
    &- \vN\cdot\left( [\frac{\vert\vu \vert^2}2 + \frac{\vert\vb \vert^2}2 +p]\vu   \right)
    + \vN \cdot ([(\vu \cdot \vb)+ q] \vb ) \\
    &+ \vu \cdot(\vN\cdot\mathbb{F}) +\vb \cdot(\vN\cdot\mathbb{G})- \mu.
\end{split}     
\end{equation}
\end{itemize}
 \end{Theorem}
 
The solutions given by Theorem \ref{weightedMHD} enjoy interesting properties as a consequence of Thorem \ref{estimates} below. \\



In the next result, we treat with discretely self-similar solutions for the (MHDG) equations. We start by remember the definition of the $\lambda$-discretely self-similarity (see \cite{CW18,PF_PG}): 
\begin{Definition} 
\begin{itemize}
\item[]
\item[$\bullet$] A vector field  $\vu_0 \in L^2_{\rm loc}(\mathbb{R}^3)$ is $\lambda$-discretely self-similar
($\vu_0$ is $\lambda$-DSS) if there exists  $\lambda>1$ such that $\lambda\vu_0(\lambda x)= \vu_0(x)$.
\item[$\bullet$] A time dependent vector field   $\vu \in L^2_{\rm loc}([0,+\infty)\times\mathbb{R}^3)$  is $\lambda$-DSS if there exists  $\lambda>1$  such that $\lambda \vu(\lambda^2 t,\lambda x)=\vu(t,x)$.
\item[$\bullet$]A forcing tensor $\mathbb{F},  \in L^2_{\rm loc}([0,+\infty)\times\mathbb{R}^3)$ is $\lambda$-DSS if   there exists $\lambda>1$  such that $\lambda^2 \mathbb{F}(\lambda^2 t,\lambda x)=\mathbb{F}(t,x)$.\\
\end{itemize}
\end{Definition} 
\begin{Theorem} \label{selfsimilar} Let $4/3<\gamma\leq 2  $ and $\lambda>1$. Let $\vu_{0}, \vb_{0}$ be $\lambda$-DSS divergence-free vector fields which belong to $ L^2_{w_\gamma}(\mathbb{R}^3)$,  and moreover, let $\mathbb{F}, \G$ be $\lambda$-DSS tensors which belong to  $ L^2_{loc} ((0,+\infty), L^2_{w_\gamma})$. Then, the (MHDG) equations
has a global weak solution $(\vu, \vb, p,q)$ such that :
 \begin{itemize} 
 \item $\vu, \vb$ is a $\lambda$-DSS vector fields.
 \item[$\bullet$] for every $0<T<+\infty$, $\vu,\vb$ belong to $L^\infty((0,T), L^2_{w_\gamma})$ and $\vN\vu, \vN \vb$ belong to $L^2((0,T),L^2_{w_\gamma})$.
 \item[$\bullet$] The map $t\in [0,+\infty)\mapsto (\vu(t), \vb(t))$ is weakly continuous from $[0,+\infty)$ to $L^2_{w_\gamma}$, and is strongly continuous at $t=0$.
 \item[$\bullet$] $(\vu,\vb,p,q)$ is suitable : it verifies the local energy inequality (\ref{energloc1}). 
\end{itemize} 
\end{Theorem}

Let us emphasize that the main contribution of this work is to establish  new a priori estimates for (MHDG)  equations (see Theorem \ref{estimates} below)  and moreover,   to show that it is simple to adapt for the (MHDG) equations the method given for the (NS) equations in \cite{PF_PG}. In this setting,   we warn that the proofs of the results in sections 3, 4 and 5 and Proposition \ref{pressurechar} keep close to their analogous in \cite{PF_PG}, but  we write them in detail for the reader understanding. 



The article is organized as follows. All our results deeply base on the study of an advection-diffusion system (AD) below and this study will be done in Section \ref{sec:advec-diff}. Then, Section \ref{sec:energ-inf} is devoted to the proof of Theorem \ref{weightedMHD}.  Finally, in Section \ref{sec:self-similar} we give a proof of Theorem \ref{selfsimilar}.

\section{The advection-diffusion problem}\label{sec:advec-diff}
From now on, we focus on the setting of the weighted Lebesgue spaces $L^p_{w_\delta}$. Let us start by recalling  their definition. For $0<\gamma$ and for all $x\in\mathbb{R}^3$  we define the weight  $w_\gamma(x)=\frac 1{(1+\vert x\vert)^\gamma}$, and then and we denote $L^p_{w_\gamma}=L^p(w_\gamma(x)\, dx)$ with $1\leq p \leq +\infty$. \\

As mentioned before, all our results base on the properties of the following advection-diffusion problem: for a time $0<T<+\infty$, let  $\vv,\vc \in L^3((0,T),L^3_{w_{3\gamma/2}})$ be  time-dependent divergence free vector-fields, then we consider the following system 
\begin{equation*}\label{AD}
(\text{AD}) \left\{ \begin{array}{ll}\vspace{2mm} 
\partial_t \vu = \Delta \vu - (\vv \cdot \nabla) \vu + (\vc \cdot \nabla) \vb - \nabla p + \vN \cdot \F,  \\ \vspace{2mm}
\partial_{t} \vb = \Delta \vb - (\vv \cdot \nabla) \vb + (\vc \cdot \nabla) \vu -\vN q + 
\vN  \cdot \G 
, \\ \vspace{2mm}
\vN \cdot  \vu=0, \,  \vN \cdot \vb =0, \\ \vspace{2mm}
\vu(0,\cdot)=\vu_0, \, \vb(0,\cdot)=\vb_0,  
\end{array}
\right.     
\end{equation*}
where $(\vu, \vb, p , q)$ are the unknowns. In the following sections, we will prove all the properties of the (AD) system that we shall need later.   
\subsection{Characterisation of the terms $p$ and $q$ and some useful results}\label{sec:Charact-p-q}
In this section we give a characterisation of the pressure  $p$ and the term $q$ (analogous to that made in \cite{PF_PG}) in the (AD) system: 
\begin{Proposition}\label{Prop-charact-p-q} Let $0\leq \gamma < \frac{5}{2} $ and $0<T<+\infty$. Let  
$\mathbb{F}(t,x)=\left(F_{i,j}(t,x)\right)_{1\leq i,j\leq 3}$  and $\mathbb{G}(t,x)=\left(G_{i,j}(t,x)\right)_{1\leq i,j\leq 3}$ be 
 tensors  such that $\mathbb{F}\in L^2((0,T), L^2_{w_\gamma})$ and $\mathbb{G}\in L^2((0,T), L^2_{w_\gamma})$.  Let $\vv,\vc \in L^3((0,T),L^3_{w_{3\gamma/2}})$ be  time-dependent divergence free vector-fields.\\
  
 Let $(\vu, \vb)$ be a solution of the following advection-diffusion problem 
\begin{equation}\label{AD2}
\left\{ \begin{array}{ll}\vspace{2mm} 
\partial_t \vu = \Delta \vu - (\vv \cdot \nabla) \vu + (\vc \cdot \nabla) \vb - \nabla \Tilde{p} + \vN \cdot \F,  \\ \vspace{2mm}
\partial_{t} \vb = \Delta \vb - (\vv \cdot \nabla) \vb + (\vc \cdot \nabla) \vu -\vN \Tilde{q} + 
\vN  \cdot \G 
, \\ \vspace{2mm}
\vN \cdot  \vu=0, \,  \vN \cdot \vb =0,
\end{array}
\right.     
\end{equation}
 
 such that $\vu, \vb \in  L^\infty((0,T), L^2_{w_\gamma})$,   $\vN\vu, \vN\vb \in L^2((0,T),L^2_{w_\gamma})$, and moreover,  $\tilde{p}$ and  $\Tilde{q}$ belongs to $\mathcal{D}'( (0,T) \times \mathbb{R} ^3 )$. \\
  
Then, the gradient terms $(\vN \Tilde{p}, \vN \Tilde{q})$  are \emph{necessarily}
related to $(\vu, \vb, \vv, \vc)$   and $\mathbb{F}$ and $\G$ through the Riesz transforms $\mathcal{R}_i =\frac{\partial_i}{\sqrt{-\Delta}}$ by the formulas
$$\vN \Tilde{p} =\vN \left(\sum_{1 \leq i,j\leq 3} \Ri_i \Ri_j(u_iv_j-b_i c_j -F_{i,j}) \right),$$
 and
  $$ \vN \Tilde{q} = \vN \left(\sum_{1 \leq i,j\leq 3} \Ri_i \Ri_j(v_i b_j-c_i u_j - G_{i,j}) \right),$$
 where,   
 \begin{equation}\label{Info1}
     \sum_{1 \leq i,j\leq 3} \Ri _i \Ri_j (u_i v_j - v_i c_j), \sum_{1 \leq i,j\leq 3} \Ri _i \Ri_j (v_i b_j - c_i u_j)  \in L^{3}((0,T),L^{6/5}_{w_{\frac{6}{5}}})
 \end{equation}
 and  
 \begin{equation}\label{Info2}
  \sum_{1 \leq i,j\leq 3} \Ri _i \Ri_j F_{i,j}, \sum_{1 \leq i,j\leq 3} \Ri _i \Ri_j G_{i,j}  \in L^{2}((0,T),L^{2}_{w_\gamma}).   
 \end{equation}
\end{Proposition}
The proof of this result deeply bases on some useful technical lemmas established in \cite{PF_PG}, Section $2$ (see also \cite{Gr08,Gr09}): 
\begin{Lemma}
\label{LemRieszt}
Let $0 \leq \delta<3$ and $1<p<+\infty$. The Riesz transforms $\mathcal{R}_i$ and the Hardy--Littlewood maximal function operator $\mathcal{M}$ are bounded on $L^p_{w_\delta}$ :
$$ \|R_jf\|_{L^p_{w_\delta}}\leq C_{p,\delta} \|f\|_{L^p_{w_\delta}} \text{ and } \|\mathcal{M}_f\|_{L^p_{w_\delta}}\leq C_{p,\delta} \|f\|_{L^p_{w_\delta}}.
$$
\end{Lemma}
\label{pressurechar}
This lemma has an important corollary which allows us to study  the convolution operator with a non increasing kernel: 
  
\begin{Lemma} 
\label{indt}
Let $0 \leq \delta<3$ and $1<p<+\infty$. If $\theta\in L^1(\mathbb{R}^3)$ is a non-negative, radial function and is radially non-increasing then for all $f\in L^p_{w_\delta}$, 
  $$ \| \theta*f\|_{L^p_{w_\delta}} \leq C_{p,\delta} \|f\|_{L^p_{w_\delta}} \|\theta\|_1.$$
\end{Lemma}
With these lemmas at hand, we are able to give a proof of Proposition \ref{Prop-charact-p-q}. \\ 

\pv We define the functions $p$ and $q$ as follows:
\begin{equation*}
 p =  \sum_{1 \leq i,j\leq 3} \Ri_i \Ri_j(u_iv_j-b_i c_j -F_{i,j}) \,\, \text{and} \,\,  q = \sum_{1 \leq i,j\leq 3} \Ri_i \Ri_j(v_i b_j-c_i u_j - G_{i,j}). 
 \end{equation*} 
 Then, by the information of the functions $(\vu, \vb, \vv, \vc, \F, \G)$ given above, using interpolation, H\"older inequalities  and the  Lemma \ref{LemRieszt} (as we have $ 0 \leq \gamma < \frac{5}{2} $) we obtain (\ref{Info1}) and (\ref{Info2}). \\
 
 We will prove now that we have $\vN(\Tilde{p}-p)=0$ and $\vN(\Tilde{q}-q)=0$. Taking the divergence operator in the equations \eqref{AD2}, as the functions $(\vu, \vb, \vv,\vc)$ are divergence-free vector fields  we obtain $ \Delta (\Tilde{p}-p)=0$  and $\Delta (\Tilde{q}-q)=0$. Then, let $\alpha \in \mathcal{D}(\mathbb{R})$ be such that $\alpha (t)= 0$  for all $|t| \geq \varepsilon$ (with $\varepsilon>0$) and moreover, let $\beta\in\mathcal{D}(\mathbb{R}^3)$. Thus, we have $ (\vN \Tilde{p} *( \alpha\otimes\beta), \vN \Tilde{q} *( \alpha\otimes\beta)) \in \mathcal{D}^{'}((\varepsilon, T-\varepsilon) \times \mathbb{R}^3).$\\
 
For $t  \in (\varepsilon, T-\varepsilon)$ fix, we  define  
$$A_{\alpha,\beta,t}=( \vN \Tilde{p} * (\alpha\otimes\beta)- \vN p * (\alpha\otimes\beta)) (t,.), $$ 
$$B_{\alpha,\beta,t}=( \vN \Tilde{q} * (\alpha\otimes\beta)- \vN q * (\alpha\otimes\beta)) (t,.).$$
Then, as $ \vN \Tilde{p}$ and $\vN \Tilde{q}$ verify the equations (\ref{AD2}) and moreover, by the properties of the convolution product,  we can write  
 \begin{equation*}
 \begin{split} 
 A_{\alpha,\beta,t}=& ( \vu *(-\partial_t\alpha\otimes\beta+\alpha\otimes\Delta\beta)+  (-\vu\otimes\vv+ \vb\otimes \vc) * (\alpha\otimes\vN\beta))(t,.)\\
 &  +\mathbb{F} * (\alpha\otimes\vN\beta))(t,.) -  ( p * (\alpha\otimes\vN\beta)) (t,.),\end{split}\end{equation*} 
 and 
  \begin{equation*}
 \begin{split} 
 B_{\alpha,\beta,t}=& ( \vb *(-\partial_t\alpha\otimes\beta+\alpha\otimes\Delta\beta)+  (-\vb\otimes\vv+ \vu\otimes \vc)* (\alpha\otimes\vN\beta))(t,.)\\
 &  +\mathbb{G} * (\alpha\otimes\vN\beta))(t,.) -  ( q * (\alpha\otimes\vN\beta)) (t,.).\end{split}\end{equation*} 
Recall that  for $\varphi\in \mathcal{D}(\mathbb{R}^3)$ we have $\vert f*\varphi\vert\leq C_\varphi \mathcal{M}_f$ and then, by  Lemma \ref{LemRieszt}, we get that a convolution with a function in $\mathcal{D}(\mathbb{R}^3)$ is a bounded operator on $L^2_{w_\gamma}$ and on $L^{6/5}_{w_{6\gamma/5}}$. Thus we have   that $A_{\alpha,\beta,t}, B_{\alpha,\beta,t} \in L^2_{w_\gamma}+ L^{6/5}_{w_{6\gamma/5}}$. Moreover, for $0<\delta$ such that $\max \{ \gamma, \frac{\gamma  + 2 }{2} \}  <\delta<5/2$ , we have $A_{\alpha,\beta,t}, A_{\alpha,\beta,t} \in L^{6/5}_{w_{6\delta/5}}$; and in particular, we have that $A_{\alpha,\beta,t}$ and $B_{\alpha,\beta,t}$ are tempered  distribution.\\

With this information, and the fact that we have  $\Delta  A_{\alpha,\beta,t}=(\alpha\otimes\beta)* \vN (\Delta(\tilde p -p))(t,.)=0$, and similarly we have $ \Delta  B_{\alpha,\beta,t}=0$,  we find that  $A_{\alpha,\beta,t}$ and $B_{\alpha,\beta,t}$  are polynomials. But,  remark that for all $1<r<+\infty$ and $0< \eta < 3$, the space $L^r_{w_\eta}$ does not contain non-trivial   polynomials and then we have  $ A_{\alpha,\beta,t}= 0$ and $ B_{\alpha,\beta,t}= 0$. Finally, we use an approximation of identity $\frac 1{\epsilon^4} \alpha	(\frac t\epsilon)\beta(\frac x\epsilon)$ to obtain that $\vN(\Tilde{p}-p)=0$ and $\vN(\Tilde{q}-q)=0$.  \Endproof{} \\

We state a Sobolev type embedding which will be very useful in the next section (see Section 2 in \cite{PF_PG}). 
 \begin{Remark}\label{sobol} For $\delta \geq 0$. Let $f\in L^2_{w_\delta}$ such that $\vN f\in L^2_{w_\delta}$ then $f\in L^6_{w_{3\delta}}$ and
  $$ \|f\|_{L^6_{w_{3\delta}}}\leq C_\delta (\|f\|_{L^2_{w_\delta}}+ \|\vN f\|_{L^2_{w_\delta}}).$$
 \end{Remark}

\subsection{A priori uniform estimates for the (AD) system} 
In order to simplify the notation, for a Banach space $X \subset \mathcal{D}'$ of vector fields endowed with a norm 
$\| \cdot \|_X $, we will write
\begin{equation*}
    \|( \vu, \vv )\|^2_X = \| \vu \|^2_X + \|\vv \|^2_X,
\phantom{space} \text{ and } \phantom{space}
    \|\vN(\vu, \vv)\|^2_X = \|\vN \vu\|^2_X + \|\vN \vv\|^2_X.
\end{equation*}

 \begin{Theorem}\label{estimates}
  Let $0 \leq \gamma\leq 2$ and $0<T<+\infty$. Let  $\vu_{0}, \vb_0 \in L^2_{w_\gamma}(\mathbb{R}^3)$ be a divergence-free vector fields and let $\mathbb{F}, \G \in L^2((0,T), L^2_{w_\gamma}) $ be  two tensors $\mathbb{F}(t,x)=\left(F_{i,j}(t,x)\right)_{1\leq i,j\leq 3}$, $\mathbb{G}(t,x)=\left(G_{i,j}(t,x)\right)_{1\leq i,j\leq 3}$.  Let $\vv,\vc \in L^3((0,T),L^3_{w_{3\gamma/2}})$ be  time-dependent divergence free vector-fields.
  
 Let $(\vu, \vb,p,q)$ be a solution of the following advection-diffusion problem 
\begin{equation}\label{AD}
(\text{AD}) \left\{ \begin{array}{ll}\vspace{2mm} 
\partial_t \vu = \Delta \vu - (\vv \cdot \nabla) \vu + (\vc \cdot \nabla) \vb - \nabla p + \vN \cdot \F,  \\ \vspace{2mm}
\partial_{t} \vb = \Delta \vb - (\vv \cdot \nabla) \vb + (\vc \cdot \nabla) \vu -\vN q + 
\vN  \cdot \G 
, \\ \vspace{2mm}
\vN \cdot  \vu=0, \,  \vN \cdot \vb =0, \\ \vspace{2mm}
\vu(0,\cdot)=\vu_0, \, \vb(0,\cdot)=\vb_0.  
\end{array}
\right.     
\end{equation}
which satisfies :
 \begin{itemize} 
 \item[$\bullet$]  $\vu, \vb$ belong to $L^\infty((0,T), L^2_{w_\gamma})$ and $\vN\vu$, $\vN \vb$ belong to $L^2((0,T),L^2_{w_\gamma})$
 \item[$\bullet$] the pressure $p$ and the term $q$ are related to $\vu$, $\vb$, $\mathbb{F}$ and $\G$ through the Riesz transforms $R_i =\frac{\partial_i}{\sqrt{-\Delta}}$ by the formulas
 $$ p =\sum_{1 \leq i,j\leq 3} \Ri_i \Ri_j(u_iv_j-b_i c_j -F_{i,j})$$
 and
  $$ q =\sum_{1 \leq i,j\leq 3} \Ri_i \Ri_j(v_i b_j-c_i u_j - G_{i,j})$$
 \item[$\bullet$] the map $t\in [0,+\infty)\mapsto (\vu(t), \vb(t ))$ is weakly continuous from $[0,+\infty)$ to $L^{2}_{w_\gamma}$, and is strongly continuous at $t=0$ :
 
 \item[$\bullet$] the solution $(\vu, \vb, p,q)$ is suitable : there exist a non-negative locally finite measure $\mu$ on $(0,+\infty)\times\mathbb{R}^3$ such that
 
\begin{equation}\label{energloc}
\begin{split}
    \partial_t(\frac {\vert\vu \vert^2 
    + |\vb|^2 }2)=&\Delta(\frac {\vert\vu \vert^2 + |\vb|^2 }2)-\vert\vN \vu \vert^2 - |\vN \vb|^2 
    - \vN\cdot\left( (\frac{\vert\vu \vert^2}2 + \frac{\vert\vb \vert^2}2)\vv \right)\\
    &-\vN \cdot(p \vu) - \vN \cdot(q \vb)+ \vN \cdot ((\vu \cdot \vb) \vc) \\
    &+ \vu \cdot(\vN\cdot\mathbb{F}) +\vb \cdot(\vN\cdot\mathbb{G})- \mu.
\end{split}     
\end{equation}

 \end{itemize}
Then we have the following controls:
\begin{itemize}
     \item[$\bullet$] If $0< \gamma \leq 2$, for almost every $a\geq 0$ (including $0$) and  for all $t\geq a$, 
    \begin{equation}\label{energcontrol}
    \begin{split}
        \Vert & ( \vu , \vb )(t) \Vert^{2}_{L^{2}_{w_\gamma}} + 2 \int_{a}^{t} (\Vert \vN ( \vu,  \vb )(s) \Vert^{2}_{L^{2}_{w_\gamma}} )ds \\
        &\leq  \Vert (\vu , \vb) (a) \Vert^{2}_{L^{2}_{w_\gamma}} -\int_a^t \int \vN(  \vert\vu\vert^2 +\vert\vb\vert^2) \cdot \vN w_\gamma\,   dx \, ds
        \\
        &+ \int_{a}^{t} \int  [(\frac{\vert\vu  \vert^2}2+ \frac{|\vb  |^2}{2} ) \vv ]  \cdot \vN w_\gamma \, dx\, ds + 2 \int_{a}^{t} \int p \vu \cdot \vN w_\gamma dx \, ds \\
        &+ 2 \int_{a}^{t} \int q \vb \cdot \vN w_\gamma dx \, ds + \int_{a}^{t} \int  [(\vu  \cdot \vb ) \vc ]  \cdot \vN w_\gamma \, dx\, ds\\
        &- \sum_{1 \leq i,j \leq 3} (\int_{a}^{t} \int   F_{i,j} (\partial_i u_j)  w_\gamma \, dx\, ds + \int_{a}^{t}\int   F_{i,j} u_j \partial_i(w_\gamma)   \, dx\, ds )\\
        &- \sum_{1 \leq i,j \leq 3} ( \int_{a}^{t}\int   G_{i,j} (\partial_i b_j)  w_\gamma \, dx\, ds + \int_{a}^{t} \int   G_{i,j} b_j \partial_i(w_\gamma)  \, dx\, ds ) ,
    \end{split}\end{equation}
    which implies in particular that the map $t \mapsto (\vu(t), \vb(t ))$  from $[0,+\infty)$ to $L^{2}_{w_\gamma}$ is stronly continuous almost everywhere and
 \begin{equation}\label{energcontrol2} 
 \begin{split}  
    \|  (\vu, \vb)(t) \|^2_{L^2_{w_\gamma}} & + \int_{a}^{t} \Vert \vN (\vu , \vb) (s) \Vert^2_{L^2_{w_\gamma}} 	 ds \\
    \leq & \|(\vu, \vb)(a) \|_{L^2_{w_\gamma}}^2  +C_\gamma \int_{a}^{t} \| ( \mathbb{F},  \mathbb{G}) (s)\|^2_{L^2_{w_\gamma}} \, ds \\
     & + C_\gamma  \int_{a}^{t}  (1+  \| (\vv, \vc )(s )\|_{L^3_{w_{3\gamma/2}}}^2 ) (\|(\vu,\vb)(s)\|_{L^2_{w_\gamma}}^2 )\, ds.
\end{split}\end{equation}

\item[$\bullet$] Si $\gamma = 0$, for almost all $a\geq 0$  (including $0$) and for all $t\geq a$, 
\begin{eqnarray*}
\begin{split}
    \Vert  (\vu, \vb)(t) \Vert^{2}_{L^2} & +  2 \int_{a}^{t} (\Vert \vN( \vu , \vb)(s) \Vert^{2}_{L^2} )	 ds \nonumber\\
     \leq&  \Vert (\vu, \vb ) (a) \Vert^{2}_{L^2} 
     \nonumber \\
     & + \sum_{1 \leq i,j \leq 3} (  \int_{a}^{t}  \int   F_{i,j}\partial_iu_{j}\    \, dx\, ds + \int_{a}^{t}  \int   G_{i,j}\partial_ib_{j}\    \, dx\, ds ) ,
\end{split}
\end{eqnarray*}
which implies of course that the map $t \mapsto (\vu(t), \vb(t ))$  from $[0,+\infty)$ to $L^{2}_{w_\gamma}$ is stronly continuous almost everywhere.
\end{itemize}
\end{Theorem}

\pv{}
We consider the case $0< \gamma \leq 2$ (the changes required for the case $\gamma = 0$ are obvious).
Let $0<t_0<t_1<T$, we take a non-decreasing function $
\alpha\in \mathcal{C}^\infty(\mathbb{R})$  equal to $0$ on $(-\infty,\frac{1}{2})$ and equal to $1$ on $(1, +\infty)$.  For $0<\eta< \min(\frac {t_0}2,T-t_1) $, let
\begin{equation}
\label{alpha}
    \alpha_{\eta,t_0,t_1}(t)=\alpha( \frac{t-t_0}\eta)-\alpha(\frac{t-t_1}\eta) .
\end{equation}
Remark that $\alpha_{\eta,t_0,t_1}$ converges almost everywhere to $\mathds 1 _{[t_0, t_1]}$ when $\eta \to 0$  and $\partial_t \alpha_{\eta,t_0,t_1}$ is the difference between two identity approximations, the first one in $t_0$ and the second one in $t_1$.\\

Consider  a non-negative function $\phi\in\mathcal{D}(\mathbb{R}^3)$ which is equal to $1$ for $\vert x\vert\leq 1$ and to $0$ for $\vert x\vert\geq 2$. We define
\begin{equation}
\label{phi}
    \phi_R(x)=\phi(\frac x R).
\end{equation}

For $\epsilon>0$, we let $w_{\gamma,\epsilon}= \frac 1{(1+\sqrt{\epsilon^2+\vert x\vert^2})^\gamma}$ (if $\gamma=0$, we let $w_{\gamma,\epsilon}=1$ ).

We have $\alpha_{\eta,t_0,t_1}(t)\phi_R(x) w_{\gamma,\epsilon}(x)\in\mathcal{D}((0,T)\times\mathbb{R}^3)$ and $\alpha_{\eta,t_0,t_1}(t)\phi_R(x) w_{\gamma,\epsilon}(x) \geq 0$. Thus, using the local energy balance (\ref{energloc}) and the fact that the measure $\mu$ verifies $\mu\geq 0$, we find

\begin{equation*} \begin{split}
  -\iint & \frac{\vert \vu \vert^2}2+ \frac{\vert \vb \vert^2}2  \partial_t\alpha_{\eta,t_0,t_1} \phi_R w_{\gamma,\epsilon} \, dx\, ds + \iint \vert\vN\vu \vert^2 + \vert\vN\vb \vert^2 \, \,  \alpha_{\eta,t_0,t_1} \phi_R w_{\gamma,\epsilon}  dx\, ds
\\ \leq &-\sum_{i=1}^3 \iint (\partial_i\vu \cdot \vu + \partial_i \vb \cdot \vb ) \,  \alpha_{\eta,t_0,t_1}  (w_{\gamma,\epsilon}\partial_i \phi_R+\phi_R \partial_i w_{\gamma,\epsilon}) \, dx\, ds
 \\
 &+ \sum_{i=1}^3 \iint  [(\frac{\vert\vu \vert^2}2+ \frac{|\vb_n |^2}{2} ) v_{ i} + p u_{i} ] \alpha_{\eta,t_0,t_1}  (w_{\gamma,\epsilon}\partial_i \phi_R+\phi_R \partial_i w_{\gamma,\epsilon}) \, dx\, ds
\\
 &+ \sum_{i=1}^3 \iint  [(\vu \cdot \vb) c_{ i}  + q b_{i} ] \alpha_{\eta,t_0,t_1}  (w_{\gamma,\epsilon}\partial_i \phi_R+\phi_R \partial_i w_{\gamma,\epsilon}) \, dx\, ds 
\\&- \sum_{1 \leq i,j \leq 3} ( \iint  F_{i,j} u_{j}  \alpha_{\eta,t_0,t_1}  (w_{\gamma,\epsilon}\partial_i \phi_R+\phi_R \partial_i w_{\gamma,\epsilon}) \, dx\, ds + \iint   F_{i,j}\partial_iu_{j}\  \alpha_{\eta,t_0,t_1} \phi_R \, dx\, ds )
\\&- \sum_{1 \leq i,j \leq 3} ( \iint  G_{i,j} b_{j}  \alpha_{\eta,t_0,t_1}  (w_{\gamma,\epsilon}\partial_i \phi_R+\phi_R \partial_i w_{\gamma,\epsilon}) \, dx\, ds + \iint   G_{i,j}\partial_ib_{j}\  \alpha_{\eta,t_0,t_1} \phi_R \, dx\, ds ) .
\end{split}\end{equation*}
As $\gamma\leq 2$, there exists $C_\gamma > 0$ which does not depend on $R>1$ nor on $\epsilon>0$
$$ \vert  w_{\gamma,\epsilon}\partial_i \phi_R\vert +\vert \phi_R \partial_iw_{\gamma,\epsilon}\vert \leq C_\gamma \frac{w_\gamma(x)}{1+\vert x\vert}\leq C_\gamma w_{3\gamma/2}(x).
$$
By interpolation we find $\vu, \vb$ belong to $L^4((0,T),L^3_{w_{3\gamma/2}})$. Also, we have $p u_i, \, q b_i \in L^1_{w_{3\gamma/2}}$ since $w_\gamma p, \, w_\gamma q \in L^2 ((0,T), L^{6/5}+L^2)$ and $w_{\gamma/2} \vu, \,  w_{\gamma/2} \vb \in L^2((0,T), L^2\cap L^6)$. Later, we will use dominated convergence using this remarks. First, we let $\eta$ go to $0$ and we find that 

 \begin{equation*} \begin{split}
  - \lim_{\eta\rightarrow 0} & \iint \frac{\vert \vu \vert^2}2+ \frac{\vert \vb \vert^2}2 \partial_t\alpha_{\eta,t_0,t_1}  \phi_R \, dx\, ds +  \int_{t_0}^{t_1}  \int \vert\vN \vu  \vert^2 + \vert\vN \vb  \vert^2 \, \,   \phi_R w_{\gamma,\epsilon}  dx\, ds \\ 
    \leq
    &-\sum_{i=1}^3 \int_{t_0}^{t_1}  \int (\partial_i\vu \cdot \vu + \partial_i \vb \cdot \vb ) \,     (w_{\gamma,\epsilon}\partial_i \phi_R+\phi_R \partial_iw_{\gamma,\epsilon}) \, dx\, ds
    \\
    &+ \sum_{i=1}^3 \int_{t_0}^{t_1}  \int  [(\frac{\vert\vu \vert^2}2+ \frac{|\vb |^2}{2} ) v_{ i} + p u_{i} ]    (w_{\gamma,\epsilon}\partial_i \phi_R+\phi_R \partial_iw_{\gamma,\epsilon}) \, dx\, ds
    \\
     &+ \sum_{i=1}^3 \int_{t_0}^{t_1}  \int  [(\vu \cdot \vb) c_{ i } + q b_{i} ]    (w_{\gamma,\epsilon}\partial_i \phi_R+\phi_R \partial_iw_{\gamma,\epsilon}) \, dx\, ds 
    \\&- \sum_{1 \leq i,j \leq 3} ( \int_{t_0}^{t_1}  \int  F_{i,j} u_{j}     (w_{\gamma,\epsilon}\partial_i \phi_R+\phi_R \partial_iw_{\gamma,\epsilon}) \, dx\, ds + \int_{t_0}^{t_1}  \int   F_{i,j}\partial_iu_{j}\    \phi_R \, dx\, ds )
    \\&- \sum_{1 \leq i,j \leq 3} ( \int_{t_0}^{t_1}  \int  G_{i,j} b_{j}    (w_{\gamma,\epsilon}\partial_i \phi_R+\phi_R \partial_iw_{\gamma,\epsilon}) \, dx\, ds + \int_{t_0}^{t_1}  \int   G_{i,j}\partial_ib_{j}\    \phi_R \, dx\, ds )
\end{split}\end{equation*}
when the limit in the left side exists. Let
 $$ A_{R,\epsilon}(t)=\int ( \vert \vu(t,x)\vert^2  + |\vb(t,x)  \vert^2 )   \phi_R(x)  w_{\gamma,\epsilon}(x) \, dx,$$
since
$$-\iint (\frac{|\vu|^2}2 + \frac{ |\vb|^2}2 ) \partial_t\alpha_{\eta,t_0,t_1}   \phi_R w_{\gamma,\epsilon} \, dx\, ds=-\frac 1 2\int \partial_t\alpha_{\eta,t_0,t_1}
A_{R,\epsilon}(s) \, ds$$
We have for all $t_0$ and $t_1$ Lebesgue points of the measurable functions $A_{R,\epsilon}$, $$ \lim_{\eta\rightarrow 0}  -\iint ( \frac{\vert\vu\vert^2}2 + \frac{\vert\vb\vert^2}2 )  \partial_t\alpha_{\eta,t_0,t_1}  \phi_R w_{\gamma,\epsilon} \, dx\, ds=\frac 1 2 (  A_{R,\epsilon}(t_1)- A_{R,\epsilon}(t_0)),$$ 

Then, by continuity, we can let $t_0$ go to $0$ and thus  replace $t_0$ by $0$ in the inequality. Moreover, if we let $t_1$ go to $t$, then by weak continuity,
we find that $$A_{R,\epsilon}(t)\leq \lim_{t_1\rightarrow t }  A_{R,\epsilon}(t_1),$$
so that we may as well replace $t_1$ by $t\in (t_0,T)$. Thus we find that for almost every $a\in (0,T)$ (including $0$) and for all $t \in (a,T)$, we have:  
\begin{equation}\label{eqlarga} 
\begin{split}
    & \frac 1 2 (  A_{R, \epsilon}(t)- A_{R, \epsilon}(a)) + \int_{a}^{t} \int \vert\vN\vu  \vert^2 + \vert\vN\vb  \vert^2 \, \,   \phi_R w_{\gamma,\epsilon} dx\, ds \\ 
    &=-\sum_{i=1}^3 \int_{a}^{t} \int (\partial_i\vu  \cdot \vu  + \partial_i \vb  \cdot \vb  ) \,    (w_{\gamma,\epsilon}\partial_i \phi_R+\phi_R \partial_iw_{\gamma,\epsilon}) \, dx\, ds \\
    &+ \sum_{i=1}^3 \int_{a}^{t} \int  [(\frac{\vert\vu  \vert^2}2+ \frac{|\vb  |^2}{2} ) v_{ i} + p  u_{i} ]   (w_{\gamma,\epsilon}\partial_i \phi_R+\phi_R \partial_iw_{\gamma,\epsilon}) \, dx\, ds \\
    &+ \sum_{i=1}^3 \int_{a}^{t} \int  [(\vu  \cdot \vb ) c_{ i } + q  b_{i} ]   (w_{\gamma,\epsilon}\partial_i \phi_R+\phi_R \partial_iw_{\gamma,\epsilon}) \, dx\, ds  \\
    &- \sum_{1 \leq i,j \leq 3} ( \int_{a}^{t} \int  F_{i,j} u_{j}    (w_{\gamma,\epsilon}\partial_i \phi_R+\phi_R \partial_iw_{\gamma,\epsilon})  \, dx\, ds - \int_{a}^{t} \int   F_{i,j}\partial_i u_{j}\   \phi_R \, dx\, ds ) \\
    &- \sum_{1 \leq i,j \leq 3} ( \int_{a}^{t} \int  G_{i,j} b_{j}    (w_{\gamma,\epsilon}\partial_i \phi_R+\phi_R w_{\gamma,\epsilon} \partial_iw_{\gamma,\epsilon}) \, dx\, ds - \int_{a}^{t} \int   G_{i,j}\partial_i b_{ j}\   \phi_R w_{\gamma,\epsilon} \, dx\, ds ) ,
\end{split}\end{equation}
Taking the limit when $R$ go to $+\infty$ and then $\epsilon$ go to $0$, by dominated convergence we obtain the energy control (\ref{energcontrol}).  We let $t$ go to $a$ in (\ref{energcontrol}), so that
 $$\limsup_{t\rightarrow 0}  \|(\vu, \vb)(t )\|_{L^2_{w_\gamma}}^2\leq  \|(\vu, \vb)(a ) \|_{L^2_{w_\gamma}}^2 .$$   Also, as $\vu$ is weakly continuous in $L^2_{w_\gamma}$,
  $$    \|(\vu, \vb)(a) \|_{L^2_{w_\gamma}}^2 \leq   \liminf_{t\rightarrow 0}  \|(\vu, \vb)(t )\|_{L^2_{w_\gamma}}^2 .$$   
  Thus $    \|(\vu, \vb)(a ) \|_{L^2_{w_\gamma}}^2 =   \lim_{t\rightarrow 0}  \|(\vu, \vb)(t)\|_{L^2_{w_\gamma}}^2 $, as we work in a Hilbert space, this fact and the weak continuity of the map $t \mapsto \vu(t) \in L^2_{w_\gamma}$ implies strongly continuity almost everywhere.\\ 

Now, to obtain \eqref{energcontrol2}, in the  energy control (\ref{energcontrol}) we have the following estimates: 
 \begin{equation*}\begin{split} \left\vert \int_0^t\int \vN\vert \vu\vert^2\cdot \vN w_\gamma\, ds\, ds\right\vert\leq&  2\gamma \int_0^t\int  \vert\vu\vert \vert\vN\vu\vert \, w_\gamma\, dx\, ds \\ \leq& \frac 1 4 \int_0^t  \|\vN\vu\|_{L^2_{w_\gamma}}^2\, ds+4\gamma^2 \int_0^t \|\vu\|_{L^2_{w_\gamma}}^2\, ds,
 \end{split}\end{equation*}
and 
$$ \left\vert \int_0^t\int \vN\vert \vb\vert^2\cdot \vN w_\gamma\, ds\, ds\right\vert\leq  \frac 1 4 \int_0^t  \|\vN\vb\|_{L^2_{w_\gamma}}^2\, ds+4\gamma^2 \int_0^t \|\vb\|_{L^2_{w_\gamma}}^2\, ds.$$
Then, for the pressure terms $p$ and $q$ we write $p=p_1+ p_2$ and $q=q_1+q_2$ where 
$$ p_1 =  \sum_{i=1}^3\sum_{j=1}^3 R_iR_j(v_iu_j - c_i b_j), \quad   p_2=-\sum_{i=1}^3\sum_{j=1}^3 R_iR_j(F_{i,j}),$$
and 
$$ q_1 =  \sum_{i=1}^3\sum_{j=1}^3 R_iR_j(v_ib_j - c_i u_j), \quad   q_2=-\sum_{i=1}^3\sum_{j=1}^3 R_iR_j(G_{i,j}),$$ 
Since $w_{6\gamma/5}\in  \mathcal{A}_{6/5}$ we have the following control
 \begin{equation*}\begin{split} & \left\vert \int_0^t\int  (\vert\vu\vert^2 \vv + |\vb|^2 \vv + ((\vu \cdot \vb) \vc) + 2p_1\vu + 2 q_1 \vb) \cdot \vN(w_\gamma)\, dx\, ds\right\vert \\
 &\leq   \gamma \int_0^t\int   (\vert\vu\vert ^2 \vert  \vv \vert + \vert\vb\vert ^2 \vert  \vv \vert + |\vu||\vb||\vc|  + 2\vert p_1\vert |\vu| + 2\vert q_1 | \vb | ) \, w_\gamma^{3/2}\, dx\, ds \\ 
 \leq &
  C_\gamma \int_0^t  \|w_\gamma^{1/2} \vu\|_6  ( \|  w_\gamma   \vert  \vv  \vert \vert\vu\vert \|_{6/5} + \| w_\gamma |\vc| \, |\vb|  \|_{6/5} )  \, ds \\
  &+ C_\gamma \int_0^t  \|w_\gamma^{1/2} \vb\|_6  ( \|  w_\gamma   \vert  \vb  \vert \vert\vv\vert \|_{6/5} + \| w_\gamma |\vc| \, |\vu|  \|_{6/5} )  \, ds \\ 
  \leq & \frac 1 4 \int_0^t  \|\vN\vu\|_{L^2_{w_\gamma}}^2\, ds + C_\gamma \int_0^t \|\vu\|_{L^2_{w_\gamma}}^2  \|\vv\|_{L^3_{w_{3\gamma/2}}}^2 + \|\vu\|_{L^2_{w_\gamma}}^2  \|\vv\|_{L^3_{w_{3\gamma/2}}}\, ds \\
   &  + C_\gamma \int_0^t \|\vb\|_{L^2_{w_\gamma}}^2  \|\vc\|_{L^3_{w_{3\gamma/2}}}^2 + \|\vu\|_{L^2_{w_\gamma}} \|\vb\|_{L^2_{w_\gamma}}
   \|\vc\|_{L^3_{w_{3\gamma/2}}}\, ds \\
  & + \frac 1 4 \int_0^t  \|\vN\vb\|_{L^2_{w_\gamma}}^2\, ds + C_\gamma \int_0^t \|\vb\|_{L^2_{w_\gamma}}^2  \|\vv\|_{L^3_{w_{3\gamma/2}}}^2 + \|\vb\|_{L^2_{w_\gamma}}^2  \|\vv\|_{L^3_{w_{3\gamma/2}}}\, ds \\
   &  + C_\gamma \int_0^t \|\vu\|_{L^2_{w_\gamma}}^2  \|\vc\|_{L^3_{w_{3\gamma/2}}}^2 + \|\vb\|_{L^2_{w_\gamma}} \|\vu\|_{L^2_{w_\gamma}}
   \|\vc\|_{L^3_{w_{3\gamma/2}}}\, ds \\
 \end{split}\end{equation*}
 and since $w_\gamma\in\mathcal{A}_2$
\begin{equation*}\label{eq06}
\begin{split}
& \left\vert   \int_{a}^{t} \int p_2  \vu \cdot \vN w_\gamma dx \, ds  + \int_{a}^{t} \int q_2  \vb \cdot \vN w_\gamma dx \, ds \right\vert  \\
\leq &  C_\gamma \int_{a}^{t} \int \vert p_2 \vert \vert \vu \vert w_\gamma \, dx \, ds  + C_\gamma \int_{a}^{t} \int \vert q_2 \vert \vert \vb \vert w_\gamma \, dx \, ds \\
\leq &  C_\gamma \int_{a}^{t}( \Vert \vu \Vert^{2}_{L^{2}_{w_\gamma}}+ \Vert p_2  \Vert^{2}_{L^{2}_{w_\gamma}} ) \, ds  + C_\gamma \int_{a}^{t} \Vert \vb \Vert^{2}_{L^{2}_{w_\gamma}}+ \Vert q_2  \Vert^{2}_{L^{2}_{w_\gamma}}  \, ds. \\
\end{split}
\end{equation*}
For the other terms, we have 
\begin{equation*}
\begin{split}
\left\vert \sum_{1 \leq i,j \leq 3} (\int_{a}^{t} \int  (F_{i,j} (\partial_i u_j)  w_\gamma \, + F_{i,j} u_j \partial_i(w_\gamma))\, dx\, ds \right\vert \leq C_\gamma \int_{a}^{t} \int \vert \F \vert (\vert \vN \vu \vert + \vert \vu \vert) w_\gamma\, dx\, ds \\
\leq \frac{1}{4} \int_{a}^{t} \Vert \vN \vu \Vert^{2}_{L^{2}_{w_\gamma}}\, ds + C_{\gamma} \int_{a}^{t} \Vert \vu \Vert^{2}_{L^{2}_{w_\gamma}}\, ds + C_{\gamma} \int_{a}^{t} \Vert \F \Vert^{2}_{L^{2}_{w_\gamma}}\, ds,
\end{split}
\end{equation*}
and 
\begin{equation*}
\begin{split}
\left\vert \sum_{1 \leq i,j \leq 3} (\int_{a}^{t} \int  G_{i,j} (\partial_i b_j)  w_\gamma \, + G_{i,j} b_j \partial_i(w_\gamma))\, dx\, ds \right\vert \leq C_\gamma \int_{a}^{t} \int \vert \F \vert (\vert \vN \vu \vert + \vert \vu \vert) w_\gamma\, dx\, ds \\
\leq \frac{1}{4} \int_{a}^{t} \Vert \vN \vb \Vert^{2}_{L^{2}_{w_\gamma}}\, ds + C_{\gamma} \int_{a}^{t} \Vert \vb \Vert^{2}_{L^{2}_{w_\gamma}}\, ds + C_{\gamma} \int_{a}^{t} \Vert \G \Vert^{2}_{L^{2}_{w_\gamma}}\, ds.
\end{split}
\end{equation*}
Hence we have found the estimate \eqref{energcontrol2} and Theorem \ref{stability} is proven.\Endproof

\section{Consequence of Gr\"onwall type inequalities and the a priori estimates.} 
 \subsection{Control for passive transportation.}
Using the Gr\"onwall inequalities, the following corollary is a direct consequence of Theorem \ref{estimates}: 
 \begin{Corollary}\label{passive} Under the assumptions of Theorem \ref{estimates}, we have
\begin{equation*}
\begin{split}
    &\sup_{0<t<T} \|(\vu,\vb) \|^2_{L^2_{w_\gamma}} + \|\vN (\vu, \vb ) \|_{L^2((0,T),L^2_{w_\gamma})} \\
     &\leq \left( \|(\vu_0,\vb_0) \|^2_{L^2_{w_\gamma}}+ C_\gamma ( \|(\mathbb{F},\mathbb{G})\|_{L^2((0,T), L^2_{w_\gamma})} )   \right) \ e^{C_\gamma (T+ T^{1/3}  \|(\vv,\vc)\|_{L^3((0,T), L^3_{w_{3\gamma /2}})}^2 )}
\end{split}
\end{equation*}
 \end{Corollary}
 Another direct consequence is the following uniqueness result for the advection-diffusion problem  (AD). 
 
  \begin{Corollary}\label{unique}. 
  Let $0\leq \gamma\leq 2$. Let $0<T<+\infty$. Let  $\vu_{0},\vb_{0} \in L^2_{w_\gamma}(\mathbb{R}^3)$ be divergence-free vector fields  and   $\mathbb{F}(t,x) =\left(F_{i,j}(t,x)\right)_{1\leq i,j\leq 3}$ and $\mathbb{G}(t,x) =\left(F_{i,j}(t,x)\right)_{1\leq i,j\leq 3}$ be  tensors such that $\mathbb{F}(t,x), \G \in L^2((0,T), L^2_{w_\gamma})$.  Let $\vv, \vc \in L^3((0,T),L^3_{w_{3\gamma/2}})$ be a time-dependent divergence free vector-fields. Assume moreover that $\vv, \vc \in L^2_t L^\infty_x(K)$ for every compact subset $K$ of  $(0,T)\times\mathbb{R}^3 $. 
 
  Let $(\vu_1, \vb_1 ,p_1,q_1)$ and $(\vu_1, \vb_1 ,p_1,q_1)$ be two solutions of the advection-diffusion problem 
\begin{equation*}
\left\{ \begin{array}{ll}\vspace{2mm} 
\partial_t \vu = \Delta \vu - (\vv \cdot \nabla) \vu + (\vc \cdot \nabla) \vb - \nabla p + \vN \cdot \F,  \\ \vspace{2mm}
\partial_{t} \vb = \Delta \vb - (\vv \cdot \nabla) \vb + (\vc \cdot \nabla) \vu -\vN q + 
\vN  \cdot \G 
, \\ \vspace{2mm}
\vN \cdot  \vu=0, \,  \vN \cdot \vb =0, \\ \vspace{2mm}
\vu(0,\cdot)=\vu_0, \, \vb(0,\cdot)=\vb_0,  
\end{array}
\right.     
\end{equation*}
which satisfies for $k=1$ or $k=2$ :
 \begin{itemize} 
 \item[$\bullet$]  $\vu_k, \vb_k$ belong to $L^\infty((0,T), L^2_{w_\gamma})$ and $\vN\vu_k$, $\vN \vb_k$ belong to $L^2((0,T),L^2_{w_\gamma})$
 
 \item[$\bullet$] the terms $p_k,q_k$ satisfy
 \begin{equation*}
     p_k =\sum_{1 \leq i,j\leq 3} \Ri_i \Ri_j(u_{k,i} v_{j} -b_{k,i} c_j -F_{i,j}),
 \end{equation*}
 and
  \begin{equation*}
     q_k =\sum_{1 \leq i,j\leq 3} \Ri_i \Ri_j (v_i b_{k,j}-c_i u_{k,j} - G_{i,j}).
 \end{equation*}

 \item[$\bullet$] the map $t\in [0,+\infty)\mapsto (\vu_k(t), \vb_k(t ))$ is weakly continuous from $[0,+\infty)$ to $L^{2}_{w_\gamma}$, and is strongly continuous at $t=0$ :

\end{itemize}

 Then $(\vu_1, \vb_1 ,p_1,q_1)=(\vu_1, \vb_1 ,p_1,q_1)$.
 
\end{Corollary}
    
\pv 
We proceed as in \cite{PF_PG} (see Corollary 5).
Let $\vw=\vu_1-\vu_2$, $\vd=\vb_1-\vb_2$, $p=p_1-p_2$ and $q = q_1- q_2$.
     Then we have
\begin{equation*}
\left\{ \begin{array}{ll}\vspace{2mm} 
\partial_t \vw = \Delta \vw - (\vv \cdot \nabla) \vw + (\vc \cdot \nabla) \vd - \nabla p ,  \\ \vspace{2mm}
\partial_{t} \vd = \Delta \vd - (\vv \cdot \nabla) \vd + (\vc \cdot \nabla) \vw -\vN q 
, \\ \vspace{2mm}
\vN \cdot  \vw=0, \,  \vN \cdot \vd =0, \\ \vspace{2mm}
\vu(0,\cdot)=0, \, \vb(0,\cdot)=0.  
\end{array}
\right.     
\end{equation*}
For all compact subset $K$ of $(0,T)\times\mathbb{R}^3$, $\vw\otimes\vv$, $ \vd\otimes\vc$, $ \vd\otimes\vv$ and $\vc\otimes\vw$ are in $L^2_t L^2_x$, and these terms belong  to $L^{3}((0,T), L^{6/5}_{w_{6\gamma/5}})$. 

We will verify that $\partial_t \vu$ and $\partial_t \vb$ are locally $L^2 H ^{-1}$. Let $\varphi, \psi\in\mathcal{D}((0,T)\times\mathbb{R}^3)$ such that $\psi=1$ on the neigborhood of the support of $\varphi$. Then
$$ \varphi p=\varphi  \Ri \otimes \Ri (\psi (\vv \otimes \vw - \vc \otimes \vd) )+\varphi  \Ri \otimes \Ri ((1-\psi )(\vv \otimes \vw - \vc \otimes \vd)).$$ 
We have that 
$$\|\varphi  \Ri \otimes \Ri (\psi (\vv \otimes \vw - \vc \otimes \vd) )\|_{L^2L^2}\leq C_{\varphi,\psi} \|\psi (\vv \otimes \vw - \vc \otimes \vd)\|_{L^2 L^2}$$ and
$$ \|\varphi  \Ri \otimes \Ri ((1-\psi )(\vv \otimes \vw - \vc \otimes \vd))\|_{L^3 L^\infty} \leq C_{\varphi,\psi} \| (\vv \otimes \vw - \vc \otimes \vd) \|_{L^3 L^{6/5}_{w_{6\gamma/5}}}
$$
with
$$ C_{\varphi,\psi}\leq C \|\varphi\|_\infty \|1-\psi\|_\infty  \sup_{x\in {\rm Supp}\, \varphi}  \left( \int_{y\in {\rm Supp }\, (1-\psi)}  \left( \frac { (1+\vert y\vert)^\gamma} {\vert x-y\vert^3}\right)^6 \right)^{1/6}<+\infty,$$
and we have analogue estimates for $\varphi q$.
Thus, we may take the scalar product of  $\partial_t \vw$ with $\vw$ and $\partial_t \vd$ with $\vd$ and find  that
\begin{equation*}
\begin{split}
    \partial_t(\frac {\vert\vw \vert^2 
    + |\vd|^2 }2)=&\Delta(\frac {\vert\vw \vert^2 + |\vd|^2 }2)-\vert\vN \vw \vert^2 - |\vN \vd|^2 
    - \vN\cdot\left( (\frac{\vert\vw \vert^2}2 + \frac{\vert\vd \vert^2}2)\vv \right)\\
    &-\vN \cdot(p \vw) - \vN \cdot(q \vd)+ \vN \cdot ((\vw \cdot \vd) \vc) \\
    &+ \vw \cdot(\vN\cdot\mathbb{F}) +\vd \cdot(\vN\cdot\mathbb{G}).
\end{split}     
\end{equation*}
 The assumptions of Theorem \ref{estimates} are satisfied then we use Corollary \ref{passive} to find that $\vw=0$ and $\vb=0$ and consequently $p=0$ and $q=0$. \Endproof

\subsection{Control for active transportation.}
We remember the following lemma (see \cite{PF_PG}) :
\begin{Lemma}\label{gronwallnl} If $\alpha$ is a non-negative bounded measurable  function on $[0,T)$ which satisfies, for two constants $A,B\geq 0$,
$$ \alpha(t)\leq A + B\int_0^t 1+\alpha(s)^3\, ds.$$ If $T_0>0$ and $T_1=\min(T,T_0, \frac 1{4B (A+BT_0)^2})$, we have, for every $t\in [0,T_1]$,  $\alpha(t)\leq \sqrt{ 2} (A+BT_0)$.
\end{Lemma}
The proof is simple, we suppose $A>0$ or $B>0$ otherwise it is obvious, let
$$\Phi(t)= A+B T_0 + B \int_0^t \alpha^3\, ds \text{ and } \Psi(t)=A+BT_0+ B \int_0^t  \Phi (s)^3 \, ds,$$ 
so that for all $t\in [0,T_1]$,  $\alpha\leq \Phi\leq \Psi$, and then
$$ \Psi'(t)= B \Phi(t)^3 \leq B \Psi(t)^3$$
so
$$ \frac 1 {\Psi(0)^2}-\frac 1{\Psi(t)^2}\leq 2Bt,$$ which let us to conclude
$$ \Psi(t)^2 \leq \frac{\Psi(0)^2}{1-2B \Psi(0)^2 t}\leq 2 \Psi(0)^2.$$
\Endproof

Now we able to prove the following tactical result.
\begin{Corollary}\label{active}  Under the hypothesis of Theorem \ref{estimates}. Assume that $(\vv, \vc)$ is controlled by $(\vu, \vb)$ in the following sense: for every $t\in (0,T)$,
$$ \| (\vv, \vc)(t)\|^2_{L^3_{w_{3\gamma/2}}}\leq C_0 \|(\vu, \vb)(t)\|^2_{L^3_{w_{3\gamma/2}}}.$$
Then there exists a constant $C_\gamma\geq 1$ such that if $T_0<T$ is such that
$$ C_\gamma  \left(1+\|(\vu_0, \vb_0 ) \|_{L^2_{w_\gamma}}^2+\int_0^{T_0} \|(\mathbb{F},\G ) \|_{L^2_{w_\gamma}}^2\, ds\right)^2\, T_0\leq 1$$ then
\begin{align*}
    \sup_{0\leq t\leq T_0} & \| (\vu, \vb)(t)\|_{L^2_{w_\gamma}}^2 + { \int_0^{T_0} \|\vN (\vu, \vb)(s) \|_{L^2_{w_\gamma}}^2\, ds }  \\
    &\leq
    C_\gamma (1 + \|(\vu_0 , \vb_0 ) \|_{L^2_{w_\gamma}}^2 +\int_0^{T_0} \|(\mathbb{F}, \G)\|_{L^2_{w_\gamma}}^2\, ds )
\end{align*}
\end{Corollary}
\pv By \eqref{energcontrol2} we can write:  
 \begin{equation*}
 \begin{split}  
    \|  (\vu, \vb)(t) \|^2_{L^2_{w_\gamma}} & + \int_{0}^{t} \Vert \vN (\vu , \vb) (s) \Vert^2_{L^2_{w_\gamma}} 	 ds \nonumber\\
    \leq & \|(\vu, \vb)(0) \|_{L^2_{w_\gamma}}^2  +C_\gamma \int_{0}^{t} \| ( \mathbb{F},  \mathbb{G}) (s)\|^2_{L^2_{w_\gamma}} \, ds \\
     & + C_\gamma  \int_{0}^{t}  (1+  \| (\vv, \vc )(s )\|_{L^3_{w_{3\gamma/2}}}^2 ) (\|(\vu,\vb)(s)\|_{L^2_{w_\gamma}}^2 )\, ds.
\end{split}\end{equation*}

Then, as we have 
 $$  \|(\vv, \vc )(s)\|_{L^3_{w_{3\gamma/2}}}^2 \leq C_0  \|(\vu, \vb)(s )\|_{L^3_{w_{3\gamma/2}}}^2 \leq C_0 C_\gamma  \| (\vu, \vb) \|_{L^2_{w_\gamma}} (\|(\vu, \vb)\|_{L^2_{w_\gamma}}+\|\vN (\vu, \vb)\|_{L^2_{w_\gamma}}),$$
 we obtain
  \begin{equation*}  \begin{split}  & \|(\vu, \vb)(t )  \|_{L^2_{w_\gamma}}^2  + \frac{1}{2}\int \| \vN (\vu, \vb) \|^2_{L^2_{w_\gamma}} \, ds  \\
  \leq &\|(\vu_0, \vb_0) \|_{L^2_{w_\gamma}}^2  +C_\gamma \int_0^t \|  (\mathbb{F}, \mathbb{G})(s)\|^2_{L^2_{w_\gamma}}\, ds  +2 C_\gamma  \int_0^t   \|(\vu, \vb)(s )\|_{L^2_{w_\gamma}}^2  +    C_0^2  \|(\vu, \vb)(s )\|_{L^2_{w_\gamma}}^6  \, ds.  \end{split}\end{equation*}  
Finally, for $t\leq T_0$ we get
   \begin{align*}
      &\|(\vu, \vb)(t)   \|_{L^2_{w_\gamma}}^2  + \frac{1}{2}\int \| \vN (\vu, \vb) \|^2_{L^2_{w_\gamma}} \, ds  \\
      &\leq \|(\vu_0, \vb_0) \|_{L^2_{w_\gamma}}^2  +C_\gamma \int_0^{T_0} \|(\mathbb{F}, \mathbb{G})\|_{L^2_{w_\gamma}}^2\, ds + C_\gamma (1+C_0^2) \int_0^t \|(\vu, \vb)(t )\|_{L^2_{w_\gamma}}^2+ \|(\vu, \vb)(t )\|_{L^2_{w_\gamma}}^6\, ds
  \end{align*}
  and then we may conclude with Lemma \ref{gronwallnl}.\Endproof
\subsection{Stability of solutions for the (AD) system}
In this section we establish the following stability result:  
\begin{Theorem}\label{stability}
  Let $0 \leq \gamma\leq 2$. Let $0<T<+\infty$. Let  $\vu_{0,n}, \vb_{0,n} \in L^2_{w_\gamma}(\mathbb{R}^3)$ be  divergence-free vector fields. Let  $\mathbb{F}_n, \G_n  \in L^2((0,T), L^2_{w_\gamma})$ be  tensors.  Let $\vv_n, \vc_n$ be   time-dependent divergence free vector-fields such that $\vv_n, \, \vc_n  \in L^3((0,T),L^3_{w_{3\gamma/2}})$.\\

Let $(\vu_n, \, \vb_n, \, p_n\, q_n ) $ be solutions of the following advection-diffusion problems 
 \begin{equation}\label{ADn}
 (AD_n) \left\{ \begin{array}{ll}\vspace{2mm} 
 \partial_t \vu_n = \Delta \vu_n - (\vv_n \cdot \nabla) \vu_n + (\vc_n \cdot \nabla) \vb_n - \nabla p_n + \vN \cdot \F_n,  \\ \vspace{2mm}
 \partial_{t} \vb_n = \Delta \vb_n - (\vv_n \cdot \nabla) \vb_n + (\vc_n \cdot \nabla) \vu_n -\vN q_n + 
 \vN  \cdot \G_n 
 , \\ \vspace{2mm}
 \vN \cdot  \vu_n=0, \, \,   \vN \cdot \vb_n =0, \\ \vspace{2mm}
 \vu_n(0,\cdot)=\vu_{0,n}, \,\, \vb_n(0,\cdot)=\vb_{0,n}.  
 \end{array}
 \right.     
 \end{equation}
  verifying the same hypothesis of Theorem \ref{estimates}. \\  
  
 
  If $(\vu_{0,n}, \, \vb_{0,n})$ is strongly convergent to $(\vu_{0,\infty}, \, \vb_{0,\infty})$ in $L^2_{w_\gamma}$,  if the sequence $(\mathbb{F}_n, \, \G_n )$ is strongly convergent to $(\mathbb{F}_\infty, \, \G_\infty)$  in $L^2((0,T), L^2_{w_\gamma})$,  and moreover,   if the sequence $(\vv_n, \, \vc_n)$ is bounded in $L^3((0,T),   L^3_{w_{3\gamma/2}})$, then there exists $\vu_\infty, \vb_\infty, \vv_\infty, \vc_\infty, p_\infty, q_\infty$   and an increasing  sequence $(n_k)_{k\in\mathbb{N}}$ with values in $\mathbb{N}$ such that
  \begin{itemize} 
 \item[$\bullet$] $(\vu_{n_k}, \vb_{n_k})$ converges *-weakly to $(\vu_\infty, \vb_\infty)$ in $L^\infty((0,T), L^2_{w_\gamma})$, $(\vN\vu_{n_k}, \vN\vb_{n_k})$ converges weakly to $(\vN\vu_\infty, \vN\vb_\infty) $ in $L^2((0,T),L^2_{w_\gamma})$.
 \item[$\bullet$] $(\vv_{n_k}, \vc_{n_k}) $ converges weakly to $(\vv_\infty, \vc_\infty)$ in $L^3((0,T), L^3_{w_{3\gamma/2}})$, $(p_{n_k},q_{n_k}) $ converges weakly to $(p_\infty, q_\infty) $ in $L^{3}((0,T),L^{6/5}_{w_{\frac {6\gamma}5}})+L^{2}((0,T),L^{2}_{w_\gamma})$.
  \item[$\bullet$] $(\vu_{n_k}, \vb_{n_k})$ converges strongly  to $(\vu_\infty, \vb_\infty)$ in  $L^2_{\rm loc}([0,T)\times\mathbb{R}^3)$ : for every  $T_0\in (0,T)$ and every $R>0$, we have
  $$\lim_{k\rightarrow +\infty} \int_0^{T_0} \int_{\vert y\vert<R} ( \vert \vu_{n_k}(s,y)-\vu_\infty(s,y)\vert^2+ \vert \vb_{n_k}(s,y)-\vb_\infty(s,y)\vert^2) \, ds\, dy=0.$$
  \end{itemize}
 
 Moreover, $(\vu_\infty, \vb_\infty, p_\infty, q_\infty)$ is  a solution of the advection-diffusion problem 
 \begin{equation}\label{ADinf}
(AD_\infty) \left\{ \begin{array}{ll}\vspace{2mm} 
\partial_t \vu_\infty = \Delta \vu_\infty - (\vv_\infty \cdot \nabla) \vu_\infty + (\vc_\infty \cdot \nabla) \vb_\infty - \nabla p_\infty + \vN \cdot \F_\infty,  \\ \vspace{2mm}
\partial_{t} \vb_\infty = \Delta \vb_\infty - (\vv_\infty \cdot \nabla) \vb_\infty + (\vc_\infty \cdot \nabla) \vu_\infty -\vN q_\infty + 
\vN  \cdot \G_\infty 
, \\ \vspace{2mm}
\vN \cdot  \vu_\infty=0, \, \,   \vN \cdot \vb_\infty =0, \\ \vspace{2mm}
\vu_\infty(0,\cdot)=\vu_{0,\infty}, \,\, \vb_\infty(0,\cdot)=\vb_{0,\infty}.  
\end{array}
\right.     
\end{equation}
 and verify the hypothesis of Theorem \ref{estimates}.
\end{Theorem}
\pv 
By Theorem \ref{estimates} and Corollary \ref{passive},  we know that $(\vu_n, \vb_n)$ is bounded in $L^\infty((0,T), L^2_{w_\gamma})$ and $(\vN\vu_n, \vN \vb_n)$ is bounded  in $L^2((0,T), L^2_{w_\gamma})$. In particular, writing
$p_n=p_{n,1}+ p_{n,2}$  with

$$ p_{n,1} =  \sum_{i=1}^3\sum_{j=1}^3 R_iR_j(v_{n,i}u_{n,j} - c_{n,i} b_{n,j}), \quad   p_2=-\sum_{i=1}^3\sum_{j=1}^3 R_iR_j(F_{n,i,j}),$$
and $q_n=q_{n,1}+q_{n,2}$ with
$$ q_{n,1} =  \sum_{i=1}^3\sum_{j=1}^3 R_iR_j(v_{n,i}b_{n,j} - c_{n,i} u_{n,j}), \quad   q_2=-\sum_{i=1}^3\sum_{j=1}^3 R_iR_j(G_{n,i,j}),$$ 
we get that  $(p_{n,1}, q_{n,1})$ is bounded in   $L^{3}((0,T),L^{6/5}_{w_{\frac {6\gamma}5}})$ and $(p_{n,2}, q_{n,2})$ is bounded in     $L^{2}((0,T),L^{2}_{w_\gamma})$. \\

Let  $\varphi\in \mathcal{D}(\mathbb{R}^3)$. We have that   $(\varphi \vu_n, \varphi \vb_n)$ are bounded in $L^2((0,T), H^1)$. Moreover,  by equations (\ref{ADn}) and by the expressions for $p_n$ and $q_n$ above, we get that 
$(\varphi\partial_t\vu_n, \varphi\partial_t\vb_n)$ are bounded in $L^2 L^2 + L^2 W^{-1,6/5}+ L^2 H^{-1}$ and then they are bounded in  $L^2((0,T), H^{-2})$. Thus, by a Rellich-Lions lemma  there exist  $(\vu_\infty, \vb_\infty )$  and an increasing  sequence $(n_k)_{k\in\mathbb{N}}$ with values in $\mathbb{N}$ such that
 $(\vu_{n_k}, \vb_{n_k}) $ converges strongly  to $(\vu_\infty, \vb_\infty)$ in  $L^2_{\rm loc}([0,T)\times\mathbb{R}^3)$ : for every  $T_0\in (0,T)$ and every $R>0$, we have
  $$\lim_{k\rightarrow +\infty} \int_0^{T_0} \int_{\vert y\vert<R} (\vert \vu_{n_k}(s,y)-\vu_\infty(s,y)\vert^2+\vert \vb_{n_k}(s,y)-\vb_\infty(s,y)\vert^2)\, dy\, ds=0.$$
As $(\vu_n, \vb_n)$ is bounded in $L^\infty((0,T), L^2_{w_\gamma})$ and $(\vN\vu_n, \vN\vu_n)$ is bounded  in $L^2((0,T), L^2_{w_\gamma})$ we have  that $(\vu_{n_k},\vb_{n_k})$ converges *-weakly to $(\vu_\infty, \vb_\infty)$ in $L^\infty((0,T), L^2_{w_\gamma})$ and we have that $(\vN\vu_{n_k}, \vN\vu_{n_k})$ converges weakly to $(\vN\vu_\infty, \vN\vb_\infty) $ in $L^2((0,T),L^2_{w_\gamma})$.\\ 

Using the Banach--Alaoglu's theorem,  there exist   $(\vv_\infty, \vc_\infty)$  such that $(\vv_{n_k}, \vc_{n_k})$ converge weakly to $(\vv_\infty, \vc_\infty)$ in $L^3((0,T), L^3_{w_{3\gamma/2}})$. In particular, we have that the  terms   $v_{n_k, i} u_{n_k, j}$, $c_{n_k, i} b_{n_k, j}$, $v_{n_k, i} b_{n_k, j}$ and $c_{n_k, i} u_{n_k, j}$ are weakly convergent in $(L^{6/5}L^{6/5})_{\rm loc}$ and thus in $\mathcal{D}'((0,T)\times \mathbb{R}^3)$. As those terms are bounded in  $L^{3}((0,T),L^{6/5}_{w_{\frac {6\gamma}5}})$, they are weakly convergent in  $L^{3}((0,T),L^{6/5}_{w_{\frac {6\gamma}5}})$. \\

Define $p_\infty= p_{\infty,1}+p_{\infty,2}$ with  $$ p_{\infty,1} =  \sum_{i=1}^3\sum_{j=1}^3 R_iR_j(v_{\infty,i}u_{\infty,j} - c_{\infty,i} b_{\infty,j}), \quad   p_2=-\sum_{i=1}^3\sum_{j=1}^3 R_iR_j(F_{\infty,i,j}),$$
and $q_\infty= q_{\infty,1}+ q_{\infty,2}$ with 
$$ q_{\infty,1} =  \sum_{i=1}^3\sum_{j=1}^3 R_iR_j(v_{\infty,i}b_{\infty,j} - c_{\infty,i} u_{\infty,j}), \quad   q_2=-\sum_{i=1}^3\sum_{j=1}^3 R_iR_j(G_{\infty,i,j}).$$  As the Riesz transforms are bounded the spaces $L^{6/5}_{w_{\frac {6\gamma}5}}$ and  $L^2_{w_\gamma}$, we find that $(p_{n_k,1}, q_{n_k,1}) $ are weakly convergent in  $L^{3}((0,T),L^{6/5}_{w_{\frac {6\gamma}5}})$  to $(p_{\infty,1}, q_{\infty,1}) $, and moreover, we find that  $(p_{n_k,2}, q_{n_k,2})$  is strongly convergent in  $L^2((0,T),L^2_{w_\gamma})$  to $(p_{\infty,2}, q_{\infty,2})$. \\
 
With those facts, we obtain that $(\vu_\infty, p_\infty, \vb_\infty, q_\infty)$ verify the following equations  in $\mathcal{D}'((0,T)\times\mathbb{R}^3)$: 
 \begin{equation*}
 \left\{ \begin{array}{ll}\vspace{2mm} 
\partial_t \vu_\infty = \Delta \vu_\infty - (\vv_\infty \cdot \nabla) \vu_\infty + (\vc_\infty \cdot \nabla) \vb_\infty - \nabla p_\infty + \vN \cdot \F_\infty,  \\ \vspace{2mm}
\partial_{t} \vb_\infty = \Delta \vb_\infty - (\vv_\infty \cdot \nabla) \vb_\infty + (\vc_\infty \cdot \nabla) \vu_\infty -\vN q_\infty + 
\vN  \cdot \G_\infty 
, \\ \vspace{2mm}
\vN \cdot  \vu_\infty=0, \, \,   \vN \cdot \vb_\infty =0.
\end{array}
\right.     
\end{equation*}
In particular,  $(\partial_t\vu_\infty, \partial_t\vb_\infty)$ belong locally to the space $L^{2}_{t} H^{-2}_{x}$, and then these functions have  representatives such that  $t\mapsto \vu_\infty(t, \cdot)$ and $t \mapsto \vb(t,\cdot)$ which are continuous from $[0,T)$ to $\mathcal{D}'(\mathbb{R}^3)$ and coincides with $\vu_\infty(0,\cdot)+ \int_0^t \partial_t \vu_\infty \, ds$ and $\vb_\infty(0,\cdot )+ \vb_\infty(0,\cdot)+ \int_0^t  \partial_t \vb_\infty \,ds$. We have necessarily $\vu_\infty(0,\cdot)=\vu_{0,\infty}$ and $\vb_\infty(0,\cdot)=\vb_{0,\infty}$ and thus $(\vu_\infty, \vb_\infty) $ is a solution of (\ref{ADinf}). 

Next,  We define 
\begin{equation*}
\begin{split}
A_{n_k} =& - \partial_t(\frac {\vert\vu_{n_k} \vert^2 
 	+ |\vb_{n_k}|^2 }2)+ \Delta(\frac {\vert\vu_{n_k}\vert^2 + |\vb_{n_k}|^2 }2) -  \vN\cdot\left( (\frac{\vert\vu_{n_k} \vert^2}2 + \frac{\vert\vb_{n_k} \vert^2}2)\vv_{n_k} \right)\\
 	& -\vN \cdot(p_{n_k} \vu_{n_k}) - \vN \cdot(q_{n_k} \vb_{n_k})+ \vN \cdot ((\vu_{n_k} \cdot \vb_{n_k}) \vc_{n_k}) \\
 &+ \vu_{n_k} \cdot(\vN\cdot\mathbb{F}_{{n_k}}) +\vb_{n_k} \cdot(\vN\cdot\mathbb{G}_{{n_k}}),
\end{split}    
\end{equation*} 
and in order to take the limit when $n_k$ go to $\infty$, we remark that by interpolation we verify  $\sqrt{w_\gamma} \vu_{n_k}$ and $\sqrt{w_\gamma}\vb_{n_k}$ are bounded and converge weakly in $L^{10/3}L^{10/3}$ and thus $ \vu_{n_k}$ and $\vb_{n_k}$ converge locally strongly in $L^3L^3$ since the locally strongly convergence in $L^2 L^2$. We also have that $p_{n_k}$ is locally bounded in $L^{3/2}L^{3/2}$, for instance if we define  $
    a_{n_k}  = \sum_{i,j} \mathcal{R }_i \mathcal{R}_j   ( \mathds{1}_{|y|< 5R}  ( v_i u_j ))
$
and $
    b_{n_k}  = - \sum_{i,j} \mathcal{R }_i \mathcal{R}_j  (  \mathds{1}_{|y|\geq 5R}  ( v_i u_j ))
$, by the continuity of $\mathcal{R}_i$ on $L^{\frac{3}{2}}$ we have
\begin{align*}
\label{cpb2}
    \int_{|x|\leq R} |a_{n_k}|^{3/2}  dx \leq C ( \int_{|x|< 5R} |  \vv_{n_k} |^{3}  dx)^{1/2} ( \int_{|x|< 5R} |  \vu_{n_k} |^{3} \, dx  )^{1/2} 
\end{align*}
and as there exist $C>0$ such that for all $|x|\leq R$ and all $|y|\geq 5R$, the kernel $K_{i,j}$ of operator $\mathcal{R}_i \mathcal{R}_j$ satisfies
$|\mathbb{K}_{i,j} (x-y)|  \leq \frac{C }{|y|^3}$ we find
\begin{align*}
   (\int_{|x|\leq R} |b_{n_k}|^{3/2}  dx )^{2/3}
    & \leq \sum_{i,j} (\int_{|x|\leq R}  ( \int_{|y|\geq 5R}  | \mathbb{K}_{i,j} (x-y) |  \, |  v_{n_k,i}(y)  u_{n_k,j}(y) | \,  \, dy )^{3/2}  dx )^{2/3} \\
   &\leq C ( \int_{|x|\leq R}   ( \int_{|y|\geq 5R}  \frac{1}{|y|^3} |  \vv_{n_k} \otimes \vu_{n_k} |  \, dy )^{3/2}  dx  )^{2/3} \\
   & \leq C R^2 \int_{|y|\geq 5R}     \frac{1}{|y|^3} |  \vv_{n_k} \otimes \vu_{n_k} |   dy  \\
   & \leq C R^2 (\int     \frac{1}{(1+|y|)^4} |  \vv_{n_k} |^2   dy )^{1/2} ( \int    \frac{1}{(1+|y|)^2} | \vu_{n_k} |^{2}   dy )^{1/2}
\end{align*}
These remarks give $A_{n_k}$ converges to $A_\infty$ in $\mathcal{D}^{'}((0,T)\times \mathbb{R}^{3})$ where 
\begin{equation*}
\begin{split}
A_{\infty} =& - \partial_t(\frac {\vert\vu_{\infty} \vert^2 
 	+ |\vb_{\infty}|^2 }2)+ \Delta(\frac {\vert\vu_{\infty}\vert^2 + |\vb_{\infty}|^2 }2) -  \vN\cdot\left( (\frac{\vert\vu_{\infty} \vert^2}2 + \frac{\vert\vb_{\infty} \vert^2}2)\vv_{\infty} \right)\\
 	&-\vN \cdot(p_{\infty} \vu_{\infty}) - \vN \cdot(q_{\infty} \vb_{\infty})+ \vN \cdot ((\vu_{\infty} \cdot \vb_{\infty}) \vc_{\infty}) \\
 &+ \vu_{\infty} \cdot(\vN\cdot\mathbb{F}_{{\infty}}) +\vb_{\infty} \cdot(\vN\cdot\mathbb{G}_{{\infty}}). 
\end{split}    
\end{equation*}
Moreover, recall by hypothesis of this theorem we have that there exist $\mu_{n_k}$ a non-negative locally finite measure on $(0,T)\times \mathbb{R}^3$ such that 

\begin{equation*}
 \begin{split}
 \partial_t(\frac {\vert\vu_{n_k} \vert^2 
 	+ |\vb_{n_k}|^2 }2)= &\Delta(\frac {\vert\vu_{n_k}\vert^2 + |\vb_{n_k}|^2 }2)-\vert\vN \vu_{n_k} \vert^2 - |\vN \vb_{n_k}|^2 \\
& - \vN\cdot\left( (\frac{\vert\vu_{n_k} \vert^2}2 + \frac{\vert\vb_{n_k} \vert^2}2)\vv_{n_k} \right)\\
& -\vN \cdot(p_{n_k} \vu_{n_k}) - \vN \cdot(q_{n_k} \vb_{n_k})+ \vN \cdot ((\vu_{n_k} \cdot \vb_{n_k}) \vc_{n_k}) \\
& + \vu_{n_k} \cdot(\vN\cdot\mathbb{F}_{{n_k}}) +\vb_{n_k} \cdot(\vN\cdot\mathbb{G}_{n_k})- \mu_{n_k}.
 \end{split}     
 \end{equation*}
Then by definition of $A_{n_k}$ we can write $ A_{n_k}= \vert\vN \vu_{n_k} \vert^2 + |\vN \vb_{n_k}|^2 + \mu_{n_k}$, and thus we have  $\displaystyle{ A_\infty= \lim_{n_k \to + \infty}  \vert\vN \vu_{n_k} \vert^2 + |\vN \vb_{n_k}|^2 + \mu_{n_k}}$.   \\

\noindent
By weak convergence, we have for a  non-negative function $\Phi\in \mathcal{D}((0,T)\times \mathbb{R}^3)$ 
 
\begin{equation*}
\begin{split}
\iint A_\infty \Phi\, dx\, ds = & \lim_{n_k\rightarrow +\infty}\iint A_{n_k} \Phi\, dx\, ds
\geq  \limsup_{n_k \to +\infty} \iint (\vert\vN \vu_{n_k} \vert^2 + |\vN \vb_{n_k}|^2) \Phi \, dx\, ds \\ 
\geq & \iint (\vert\vN \vu_{\infty} \vert^2 + |\vN \vb_{\infty}|^2) \Phi \, dx\, ds.
\end{split}    
\end{equation*} Thus, there   exists a non-negative locally finite measure $\mu_\infty$ on $(0,T)\times\mathbb{R}^3$ such that $A_\infty=(\vert\vN \vu_\infty\vert^2+\vert\vN \vb_\infty\vert^2) +\mu_\infty$, and then we have 
 \begin{equation*}
 \begin{split}
 \partial_t(\frac {\vert\vu_\infty \vert^2 
 	+ |\vb_\infty|^2 }2)= &\Delta(\frac {\vert\vu_\infty\vert^2 + |\vb_\infty|^2 }2)-\vert\vN \vu_\infty \vert^2 - |\vN \vb_\infty|^2 \\
& - \vN\cdot\left( (\frac{\vert\vu_\infty \vert^2}2 + \frac{\vert\vb_\infty \vert^2}2)\vv_\infty \right)\\
& -\vN \cdot(p_\infty \vu_\infty) - \vN \cdot(q_\infty \vb_\infty)+ \vN \cdot ((\vu_\infty \cdot \vb_\infty) \vc_\infty) \\
& + \vu_\infty \cdot(\vN\cdot\mathbb{F}_{\infty}) +\vb_\infty \cdot(\vN\cdot\mathbb{G}_{_\infty})- \mu_\infty.
 \end{split}     
 \end{equation*}
As in \eqref{eqlarga} with the functions $(\vu_{n_k}, p_{n_k}, \vb_{n_k}, q_{n_k})$ and with $a=0$, and moreover, taking the limsup when $n_k \to +\infty$ we have 

\begin{equation*}
\begin{split}
  & \limsup_{n_k \to +\infty} \left( \int (\frac{\vert \vu_{n_k}(t,x)\vert^2 }{2}+ \frac{\vert \vb_{n_k}(t,x)\vert^2 }{2})\phi_R w_{\gamma, \varepsilon} \, dx     + \int_{0}^{t} \int \vert\vN\vu_{n_k}  \vert^2 + \vert\vN\vb_{n_k}  \vert^2 \, \,   \phi_R w_{\gamma,\epsilon} dx\, ds \right) 
\\ &\leq \int (\frac{\vert \vu_{0,\infty}(x)\vert^2 }{2}+ \frac{\vert \vb_{0,\infty}(x)\vert^2 }{2})\phi_R w_{\gamma, \varepsilon} \, dx \\
& -\sum_{i=1}^3 \int_{0}^{t} \int (\partial_i\vu_{\infty}  \cdot \vu_{\infty}  + \partial_i \vb_{\infty}  \cdot \vb_{\infty}) \,    (w_{\gamma,\epsilon}\partial_i \phi_R+\phi_R \partial_iw_{\gamma,\epsilon}) \, dx\, ds
 \\
 &+ \sum_{i=1}^3 \int_{0}^{t} \int  [(\frac{\vert\vu_{\infty}  \vert^2}2+ \frac{|\vb_{\infty}  |^2}{2} ) v_{\infty, i} + p_{\infty}  u_{\infty,i} ]   (w_{\gamma,\epsilon}\partial_i \phi_R+\phi_R \partial_iw_{\gamma,\epsilon}) \, dx\, ds
\\
 &+ \sum_{i=1}^3 \int_{0}^{t} \int  [(\vu_{\infty}  \cdot \vb_{\infty} ) c_{\infty,i} + q_{\infty}  b_{\infty,i} ]   (w_{\gamma,\epsilon}\partial_i \phi_R+\phi_R \partial_iw_{\gamma,\epsilon}) \, dx\, ds 
\\&- \sum_{1 \leq i,j \leq 3} ( \int_{0}^{t} \int  F_{\infty,i,j} u_{\infty,j}    (w_{\gamma,\epsilon}\partial_i \phi_R+\phi_R \partial_iw_{\gamma,\epsilon})  \, dx\, ds - \int_{0}^{t} \int   F_{\infty,i,j}\partial_i u_{\infty,j}\   \phi_R \, dx\, ds )
\\&- \sum_{1 \leq i,j \leq 3} ( \int_{0}^{t} \int  G_{\infty,i,j} b_{\infty,j}    (w_{\gamma,\epsilon}\partial_i \phi_R+\phi_R w_{\gamma,\epsilon} \partial_i w_{\gamma,\epsilon}) \, dx\, ds - \int_{0}^{t} \int   G_{\infty,i,j}\partial_i b_{\infty, j}\,\phi_R w_{\gamma,\epsilon} \, dx\, ds).
\end{split}\end{equation*}
Now, recall that we have $\displaystyle{\vu_{n_k}= 
\vu_{0, n_k}+ \int_0^t \partial_t \vu_{n_k}\, ds}$ and $\displaystyle{\vb_{n_k}= 
\vu_{0, n_k}+ \int_0^t \partial_t \vb_{n_k}\, ds}$ and then, for all $t \in (0,T)$ we have that $(\vu_{n_k}(t, \cdot), \vb_{n_k}(t,\cdot))$ converge to $(\vu_{\infty}(t, \cdot), \vb_{\infty}(t,\cdot))$ in $\mathcal{D}'(\mathbb{R}^3)$. Moreover, as $(\vu_{n_k}(t,\cdot), \vb_{n_k}(t,\cdot))$ are bounded in $L^2_{w_\gamma}(\mathbb{R}^3)$ we get that $(\vu_{n_k}(t,\cdot), \vb_{n_k}(t,\cdot))$ converge to $(\vu_{\infty}(t,\cdot), \vb_{\infty}(t,\cdot))$ in $L^2_{\rm loc}(\mathbb{R}^3)$. 
Thus, we can write 
\begin{equation*}
\begin{split}
 \int (\frac{\vert \vu_\infty(t,x)\vert^2}2+\frac{\vert \vb_\infty(t,x)\vert^2}2 ) \phi_R w_{\gamma,\epsilon}\, dx 
 \leq  
\limsup_{n_k \to +\infty}  \int (\frac{\vert \vu_{n_k}(t,x)\vert^2 }{2}+ \frac{\vert \vb_{n_k}(t,x)\vert^2 }{2})\phi_R w_{\gamma, \varepsilon} \, dx. 
\end{split}    
\end{equation*}
On the other hand, as $(\vN\vu_{n_k},\vN\vb_{n_k})$ are weakly convergent to $(\vN\vu_{\infty},\vN\vb_{\infty})$ in $L^{2}_{t}L^2_{w_\gamma}$, we have
\begin{equation*}
\begin{split}
&\int_0^t \int (\frac{\vert\vN \vu_{\infty}(s,x)\vert^2}{2}+\frac{\vert\vN \vu_{\infty}(s,x)\vert^2}{2}) \phi_R w_{\gamma, \varepsilon}\, dx\, ds \\
\leq & \limsup_{n_k \to +\infty} \int_{0}^{t} \int \vert\vN\vu_{n_k}  \vert^2 + \vert\vN\vb_{n_k}  \vert^2 \, \,   \phi_R w_{\gamma,\epsilon} dx\, ds. 
\end{split}    
\end{equation*}
Thus, taking the limit when $R \to 0$ and when $\varepsilon \to 0 $, for every $t\in (0,T)$ we get:

\begin{equation*}
\begin{split}
     \Vert & ( \vu_{\infty} , \vb_{\infty} )(t) \Vert^{2}_{L^{2}_{w_\gamma}} + 2 \int_{0}^{t} (\Vert \vN ( \vu_{\infty},  \vb_{\infty} )(s) \Vert^{2}_{L^{2}_{w_\gamma}} )ds \\
     \leq&  \Vert (\vu_{0,\infty} , \vb_{0,\infty}) \Vert^{2}_{L^{2}_{w_\gamma}} -\int_0^t \int (  \vN\vert\vu_{\infty}\vert^2 + \vN\vert\vb_{\infty}\vert^2) \cdot \vN w_\gamma\,   dx \, ds
    \\
     &+ \int_{0}^{t} \int  [(\frac{\vert\vu_{\infty}  \vert^2}2+ \frac{|\vb_{\infty} |^2}{2} ) \vv ]  \cdot \vN w_\gamma \, dx\, ds + 2 \int_{0}^{t} \int p_{\infty} \vu_{\infty} \cdot \vN w_\gamma dx \, ds \\
     &+ 2 \int_{0}^{t} \int q_{\infty} \vb_{\infty} \cdot \vN w_\gamma dx \, ds + \int_{0}^{t} \int  [(\vu_{\infty}  \cdot \vb_{\infty} ) \vc_{\infty} ]  \cdot \vN w_\gamma \, dx\, ds\\
     &- \sum_{1 \leq i,j \leq 3} (\int_{0}^{t} \int   F_{\infty,i,j} (\partial_i u_{\infty,j})  w_\gamma \, dx\, ds + \int_{0}^{t}\int   F_{\infty,i,j} u_{\infty,i} \partial_j(w_\gamma) \cdot \vN w_\gamma  \, dx\, ds )\\
     &- \sum_{1 \leq i,j \leq 3} ( \int_{0}^{t}\int   G_{\infty,i,j} (\partial_i b_{\infty,j})  w_\gamma \, dx\, ds + \int_{0}^{t} \int   G_{\infty,i,j} b_{\infty,i} \partial_j(w_\gamma)  \, dx\, ds ). 
\end{split}\end{equation*}
In this estimate we take now the limsup when $t\to 0$ 
in order to find that
$$ \lim_{t\rightarrow 0}    \Vert  ( \vu_{\infty} , \vb_{\infty} )(t) \Vert^{2}_{L^{2}_{w_\gamma}} =   \Vert  ( \vu_{0,\infty} , \vb_{0,\infty} ) \Vert^{2}_{L^{2}_{w_\gamma}}.$$ 
which implies strongly convergence of the solution to the initial data (since we have weak convergence and convergence of the norms in a Hilbert space). The proof is finished.
\Endproof\\

\begin{Remark}
\label{Stabilitynln}
We remark that non linear versions of this stability theorem emerge from the same proof if we take $\vu_{n}=\vv_{n}$ and $\vc_{n}=\vb_{n}$, in which case we obtain $\vu_{\infty}=\vv_{\infty}$ and $\vc_{\infty}=\vb_{\infty}$. We consider two cases.
\begin{itemize}
    \item if $\vu_{n}=\vv_{n}$ and $\vc_{n}=\vb_{n}$, the same proof give a solution on $(0,T_0)$, where $T_0<T$ using Theorem \ref{estimates} and Corollary \ref{active}.
    \item if $\vu_{n}=\vv_{n}$ and $\vc_{n}=\vb_{n}$, and we suppose that $(\vu_n, \vb_n)$ is bounded in $L^\infty((0,T), L^2_{w_\gamma})$ and $(\vN\vu_n, \vN \vb_n)$ is bounded  in $L^2((0,T), L^2_{w_\gamma})$, the same proof give a solution on $(0,T)$. We will use this case in the end of the proof of Theorem \ref{weightedMHD}.
\end{itemize}
\end{Remark}

\section{Global solutions for 3D MHD equations}\label{sec:energ-inf}

\subsection{Proof of Theorem \ref{weightedMHD}}
 Initially, we proof the local in time  existence of solutions. 
\subsubsection{Local existence}\label{sec:local-existence}
Let $\phi \in \mathcal{D}(\mathbb{R}^3)$ be a non-negative function such that $\phi(x)=1$ for $\vert x \vert <1$ and $\phi(x)=0$ for $\vert x \vert \geq 2$. For $R>0$, we define the cut-off function  $\phi_R(x)= \phi(\frac{x}{R})$. Then, for the initial $(\vu_0, \vb_0) \in L^2_{w_\gamma}(\mathbb{R}^3)$ we define  $(\vu_{0,R}, \vb_{0,R})=(\P(\phi_R \vu_0), \P(\phi_R \vb_0))\in L^{2}(\mathbb{R}^3) $ which are divergence-free  vector fields. Moreover, for the tensors $\F, \G \in L^2((0,T), L^2_{w_\gamma})$ we define $(\F_{R}, \G_R)= (\phi_R \F, \phi_R \G) \in L^2((0,T), L^2)$. \\

Then, by \ref{Prop-approx-leray}, there exist $ \vu_{R, \epsilon}, \vb_{R, \epsilon} , p_{R, \epsilon} , q_{R, \epsilon}$ solving 
\begin{equation*}
\left\{ \begin{array}{ll}\vspace{2mm} 
\partial_t \vu_{R, \epsilon} = \Delta \vu_{R, \epsilon} - ((\vu_{R,\epsilon}*\theta_\epsilon) \cdot \nabla) \vu_{R, \epsilon} + ((\vb_{R,\epsilon}*\theta_\epsilon) \cdot \nabla) \vb_{R, \epsilon} - \nabla p_{R, \epsilon} + \vN \cdot \F_R,  \\ \vspace{2mm}
\partial_{t} \vb_{R, \epsilon} = \Delta \vb_{R, \epsilon} - ((\vu_{R,\epsilon}*\theta_\epsilon) \cdot \nabla) \vb_{R, \epsilon} + ((\vb_{R,\epsilon}*\theta_\epsilon) \cdot \nabla) \vu_{R, \epsilon} -\nabla q_{R, \epsilon} + \vN  \cdot \G_R 
, \\ \vspace{2mm}
\vN \cdot  \vu_{R, \epsilon}=0, \,  \vN \cdot \vb_{R, \epsilon} =0, \\ \vspace{2mm}
\vu_{R, \epsilon}(0,\cdot)=\vu_{0,R}, \, \vb_{R, \epsilon}(0,\cdot)=\vb_{0,R}.  
\end{array}
\right.     
\end{equation*}
such that $ (\vu_{R, \epsilon}, \vb_{R, \epsilon}) \in \mathcal{C} ([0, T),L^{2}(\Rt))\cap L^{2}([0, T), \dot{H}^{1}(\Rt))$ and $(p_{R, \epsilon}, q_{R, \epsilon})  \in L^{4}((0,T),L^{6/5}(\Rt)) + L^{2}((0,T),L^{2}(\Rt))$, for every $0<T<+\infty$, and satisfying the energy equality \eqref{eneqpo}.

Now, we must study the convergence of the solution $(\vu_{R, \epsilon}, \vb_{R, \epsilon}, p_{R, \epsilon}, q_{R, \epsilon})$ when we let $R\to +\infty$ and $\epsilon \to 0$ and for this we will use the Theorem \ref{stability}, which was proven in the setting of the advection-diffusion problem  (\ref{ADn}). Thus, the first thing to do is to set   $(\vv_{R,\epsilon}, \vc_{R,\epsilon}) =(\vu_{R,\epsilon}*\theta_\epsilon,\vb_{R,\epsilon}*\theta_\epsilon)$ in (\ref{ADn}), and then, we will prove that $(\vv_{R,\epsilon}, \vc_{R,\epsilon})$ are uniform bounded in $L^{3}((0,T_0), L^{3}_{3 \gamma / 2})$ for a time $T_0>0$ small enough. \\

For a time $0<T<+\infty$, by Lemma \ref{LemRieszt} we have  
\begin{equation*}
\begin{split}
&  \|(\vv_{R, \epsilon}, \vc_{R, \epsilon} )\|_{L^3((0,T), L^3_{w_{3\gamma/2}} ) } \leq \| (\mathcal{M}_{\vu_{R,\epsilon}}, \mathcal{M}_{\vb_{R, \epsilon}} )\|_{L^3((0,T), L^3_{w_{3\gamma/2}} ) }\\
\leq & C_\gamma \| (\vu_{R,\epsilon}, \vb_{R,\epsilon} )\|_{L^3((0,T), L^3_{w_{3\gamma/2}} )}.  
\end{split}     
\end{equation*} 
Then, by interpolation and by Remark \ref{sobol} we can write 
\begin{equation*}
\begin{split}
\| (\vu_{R,\epsilon}, \vb_{R,\epsilon} )\|_{L^3((0,T), L^3_{w_{3\gamma/2}} )} \leq & C_\gamma T^{1/12}\left( (1+\sqrt{T}) \Vert (\vu_{R, \epsilon}, \vb_{R, \epsilon}) \Vert_{L^{2}((0,T), L^{2}_{w_\gamma})}  \right) \\
& + C_\gamma T^{1/12} \left( (1+\sqrt{T}) \Vert ( \vN\vu_{R, \epsilon}, \vN\vb_{R, \epsilon}) \Vert_{L^{2}((0,T), L^{2}_{w_\gamma})}  \right).  
\end{split}
\end{equation*}
At this point, remark  that $(\vu_{R, \epsilon}, \vb_{R, \epsilon}, p_{R, \epsilon}, q_{R, \epsilon})$  satisfy  the assumptions of Theorem  \ref{estimates} and then we can apply Corollary \ref{active}. Thus, for a time $T_0>0$ such that
$$ C_\gamma  \left(1+\|(\vu_{0,R}, \vb_{0,R}) \|_{L^2_{w_\gamma}}^2+\int_0^{T_0} \|(\F_R ,  \G_R) \|_{L^2_{w_\gamma}}^2\, ds\right)^2\, T_0\leq 1,$$  we have the estimates
\begin{align*}
\begin{split}
\sup_{0\leq t\leq T_0} & \| (\vu_{R,\epsilon}, \vb_{R,\epsilon})(t)\|_{L^2_{w_\gamma}}^2 + { \int_0^{T_0} \|\vN (\vu_{R,\epsilon}, \vb_{R,\epsilon})(s) \|_{L^2_{w_\gamma}}^2\, ds } \\
&\leq C_\gamma (1 + \|(\vu_{0,R}, \vb_{0,R}) \|_{L^2_{w_\gamma}}^2 +\int_0^{T_0} \|(\F_R ,  \G_R) \|_{L^2_{w_\gamma}}^2\, ds).   
\end{split}    
\end{align*}
Moreover, we have that
$$  \|(\vu_{0,R}, \vb_{0,R}) \|_{L^2_{w_\gamma}}\leq C_\gamma  \|(\vu_0, \vb_0)\|_{L^2_{w_\gamma}}  \text{ and }  \|(\F_R ,  \G_R) \|_{L^2_{w_\gamma}} \leq \|( \mathbb{F}, \G) \|_{L^2_{w_\gamma}}.$$
 and thus, by  the estimates above we find that $(\vv_{R,\epsilon}, \vc_{R,\epsilon})$ are uniform bounded in $L^3((0,T_0),   L^3_{w_{3\gamma/2}})$.\\
  
Now, we are able to apply the Theorem \ref{stability}. Let  $R_n $ be a sequence converging to $+\infty$ and $\epsilon_n$ a sequence converging to $0$  and let us denote   $(\vu_{0,n},\vb_{0,n})=(\vu_{0,R_n}, \vb_{0,R_n} )$, $ (  \mathbb{F}_n, \G_n ) = (\mathbb{F}_{R_n}, \mathbb{G}_{R_n} )$, $(\vv_n,\vc_n) =(\vv_{R_n,\epsilon_n}, \vv_{R_n,\epsilon_n})$ and $(\vu_n, \vb_n)=(\vu_{R_n,\epsilon_n},\vb_{R_n,\epsilon_n})$. As $(\vu_{0,n},\vb_{0,n})$ is strongly convergent to $(\vu_{0}, \vb_0)$ in $L^2_{w_\gamma}$,   $(\mathbb{F}_n, \G_n)$ is strongly convergent to $(\mathbb{F}, \G)$  in $L^2((0,T_0), L^2_{w_\gamma})$,  and moreover, as  $(\vv_n, \vc_n)$ is uniform  bounded in $L^3((0,T_0),L^3_{w_{3\gamma/2}})$, by Theorem \ref{stability}  there exist $(\vu, \vb, \vv, \vc, p,q)$    and there exists  an increasing  sequence $(n_k)_{k\in\mathbb{N}}$ with values in $\mathbb{N}$ such that:
\begin{itemize} 
 \item[$\bullet$] $(\vu_{n_k}, \vb_{n_k} )$ converges *-weakly to $(\vu, \vb)$ in $L^\infty((0,T_0), L^2_{w_\gamma})$, $ (\vN \vu_{n_k}, \vN \vb_{n_k} )$ converges weakly to $ (\vN \vu, \vN\vb)$ in $L^2((0,T_0),L^2_{w_\gamma})$.
 \item[$\bullet$] $(\vv_{n_k}, \vc_{n_k} )$ converges weakly to $(\vv, \vc)$ in $L^3((0,T_0), L^3_{w_{3\gamma/2}})$. Moreover,  $p_{n_k}$ converges weakly to $p$ in $L^{3}((0,T_0),L^{6/5}_{w_{\frac {6\gamma}5}})+L^{2}((0,T_0),L^{2}_{w_\gamma})$ and similarly for $q_{n_k}$.
  \item[$\bullet$] $(\vu_{n_k}, \vb_{n_k} )$ converges strongly  to $(\vu, \vb)$ in  $L^2_{\rm loc}([0,T_0)\times\mathbb{R}^3)$.  \end{itemize}

 Moreover, $(\vu, \vb, \vv, \vc, p, q)$ is  a solution of the   advection-diffusion problem 
 \begin{equation*}
 \left\{ \begin{array}{ll}\vspace{2mm} 
\partial_t \vu = \Delta \vu - (\vv \cdot \nabla) \vu + (\vc \cdot \nabla) \vb - \nabla p + \vN \cdot \F,  \\ \vspace{2mm}
\partial_{t} \vb = \Delta \vb - (\vv \cdot \nabla) \vb + (\vc \cdot \nabla) \vu -\vN q + 
\vN  \cdot \G 
, \\ \vspace{2mm}
\vN \cdot  \vu=0, \,  \vN \cdot \vb =0, \\ \vspace{2mm}
\vu(0,\cdot)=\vu_0, \, \vb(0,\cdot)=\vb_0.  
\end{array}
\right.     
\end{equation*}
 and is
  such that :
 \begin{itemize}  
 \item[$\bullet$] the map $t\in [0,T_0)\mapsto (\vu(t), \vb(t))$ is weakly continuous from $[0,T_0)$ to $L^2_{w_\gamma}$, and is strongly continuous at $t=0$
 
 \item[$\bullet$]  there exists a non-negative locally finite measure $\mu$ on $(0,T_0)\times\mathbb{R}^3$ such that
\begin{equation*}
\begin{split}
    \partial_t(\frac {\vert\vu \vert^2 
    + |\vb|^2 }2)=&\Delta(\frac {\vert\vu \vert^2 + |\vb|^2 }2)-\vert\vN \vu \vert^2 - |\vN \vb|^2 
    - \vN\cdot\left( (\frac{\vert\vu \vert^2}2 + \frac{\vert\vb \vert^2}2)\vv \right)\\
    &-\vN \cdot(p \vu) - \vN \cdot(q \vb)+ \vN \cdot ((\vu \cdot \vb) \vc) \\
    &+ \vu \cdot(\vN\cdot\mathbb{F}) +\vb \cdot(\vN\cdot\mathbb{G})- \mu.
\end{split}     
\end{equation*}
 \end{itemize} 
Finally we must check that  $\vv=\vu$ and $\vc=\vb$. As we have  $\vv_n=\theta_{\epsilon_n}*(\vu_n-\vu)+ \theta_{\epsilon_n}*\vu$, and $\vc_n=\theta_{\epsilon_n}*(\vb_n-\vb)+ \theta_{\epsilon_n}*\vb$  then we get that  $(\vv_{n_k},\vc_{n_k})$ are strongly   convergent to $(\vu, \vb)$ in   $L^3_{\rm loc}([0,T_0)\times\mathbb{R}^3)$, hence we have  $\vv=\vu$ and $\vc=\vb$. Thus, $(\vu, \vb, p ,q)$  is a solution of the (MHDG) equations  on $(0,T_0)$.

\subsubsection{Global existence}
 Let $\lambda>1$. For $n\in\mathbb{N}$ we consider the (MHDG) problem with the initial data $(\tilde \vu_{0,n}, \tilde \vb_{0,n})=(\lambda^n \vu_0(\lambda^n \cdot), \vb_0(\lambda^n \cdot))$ and with the   tensors $\mathbb{F}_n=\lambda^{2n}\mathbb{F}(\lambda^{2n} \cdot , \lambda^n \cdot )$ and $\mathbb{G}_n=\lambda^{2n}\mathbb{G}(\lambda^{2n} \cdot,  \lambda^n \cdot )$ and then, by Section \ref{sec:local-existence} we have a solution a local in time $(\tilde \vu_n, \tilde \vb_n,  \tilde p_n, \tilde q_n )$  on the interval of time  $(0,T_n)$, where the time $T_n>0$ is such that  
 \begin{equation}
 \label{eq1ss}
     C_\gamma \left(1+\|( \tilde \vu_{0,n}, \tilde \vb_{0,n}) \|_{L^2_{w_\gamma}}^2+\int_0^{+\infty} \|(\mathbb{F}_n, \G_n)\|_{L^2_{w_\gamma}}^2\, ds\right)^2\, T_n =  1
 \end{equation}

Moreover,  using the scaling of the  (MHDG) equations, which is the same well-know scaling of the Navier-Stokes equations, we can write  
$$(\tilde \vu_n, \tilde \vb_n)= (\lambda^{n}{\vu}_n(\lambda^{2n}t\cdot,\lambda^n\cdot ),\lambda^{n}{\vb}_n(\lambda^{2n}t\cdot,\lambda^n\cdot )),$$ where $({\vu}_n, {\vb}_n)$ is a solution of the (MHDG) equations on the interval of time $(0, \lambda^{2n} T_n)$ and arising from the data $(\vu_0, \vb_0, \F, \G)$. \\

By Lemma $10$ in \cite{PF_PG} we have  $\lim_{n\rightarrow +\infty} \lambda^{2n} T_n=+\infty$ and  then, for a time $T>0$ there exist $n_T \in \mathbb{N}$ such that for all $n>n_T$ we have $  \lambda^{2n}T_n>T$.  From the solution $({\vu}_n, {\vb}_n)$ on $(0,T)$ given above,  for all $n>n_T$ we define the functions 
$$(\Tilde{\Tilde{\vu}}_n, \Tilde{\Tilde{\vb}}_n)= (\lambda^{n_T}{\vu}_n(\lambda^{2n_T}t\cdot,\lambda^{n_T}\cdot ),\lambda^{n_T}{\vb}_n(\lambda^{2n_T}t\cdot,\lambda^{n_T}\cdot )),$$
which are solutions of the (MHDG) equations on $(0,\lambda^{-2n_T} T)$ with initial data $(\tilde \vu_{0,n_T},\tilde \vb_{0,n_T} )$ and forcing tensors $\mathbb{F}_{n_T},\mathbb{G}_{n_T}$, moreover, as  $\lambda^{-2n_T}T\leq T_{n_T}$,  \eqref{eq1ss} implies
  $$   C_\gamma \left(1+\|(\tilde \vu_{0,n_T},\tilde \vb_{0,n_T} )\|_{L^2_{w_\gamma}}^2+\int_0^{+\infty} \|(\mathbb{F}_{n_T},\mathbb{G}_{n_T})\|_{L^2_{w_\gamma}}^2\, ds\right)^2\, \lambda^{-2n_T} T \leq  1.$$  
and thus Corollary \ref{active} gives
\begin{align*}
    \sup_{0\leq t\leq   \lambda^{-2n_T} T} & \|\ (\tilde {\tilde \vu}_n, \tilde {\tilde \vb}_n)(t)\|_{L^2_{w_\gamma}}^2 + { \int_0^{\lambda^{-2n_T}T} \|\vN(\tilde {\tilde \vu}_n, \tilde {\tilde \vb}_n)\|_{L^2_{w_\gamma}}^2\, ds } \\
    & \leq
    C_\gamma (1+  \|(\tilde \vu_{0,n_T},\tilde \vb_{0,n_T} )\|_{L^2_{w_\gamma}}^2 +\int_0^{\lambda^{-2n_T}T} \|(\mathbb{F}_{n_T},\mathbb{G}_{n_T})\|_{L^2_{w_\gamma}}^2\, ds 
)
\end{align*}

This fact permits to find uniform controls on $(0, T)$ for $(\vu_n, \vb_n)$, with $n > n_T$, since
\begin{equation*}
\begin{split}
\lambda^{n_T(\gamma-1)} \|({\vu}_n, &{\vb}_n)(\lambda^{2n_T}t,.)\|_{L^2_{w_\gamma}}^2 \\
\leq  & \int \vert({\vu}_n(\lambda^{2n_T}t,x)\vert^2+ {\vb}_n(\lambda^{2n_T}t,x)\vert^2) \lambda^{n_T(\gamma-1)} \frac{(1+\vert x\vert)^\gamma}{(\lambda^{n_T}+\vert x\vert)^\gamma}  w_\gamma(x)\, dx \\
\leq &  \| (\Tilde{\Tilde{\vu}}_n,\Tilde{\Tilde{\vb}}_n)(t,\cdot)\|_{L^2_{w_\gamma}}^2,
\end{split}    
\end{equation*}
and 
\begin{equation*}
\begin{split}
 &\lambda^{n_T(\gamma-1)}  \int_0^T \|( \vN {\vu}_n, \vN {\vb}_n )(s,\cdot)\|_{L^2_{‘w_\gamma}}^2\, ds \\
\leq &  \int_0^T\int ( \vert\vN{\vu}_n(s,x)\vert^2+ \vert\vN{\vb}_n(s,x)\vert^2 )\lambda^{n_T(\gamma-1)} \frac{(1+\vert x\vert)^\gamma}{(\lambda^{n_T}+\vert x\vert)^\gamma}  w_\gamma(x)\, dx \, ds  \\
\leq & \int_0^{\lambda^{-2n_T}T} \|(\vN \Tilde{\Tilde{\vu}}_{n}, \vN \Tilde{\Tilde{\vb}}_{n} )(s,\cdot)\|_{L^2_{w_\gamma}}^2\, ds. 
\end{split}    
\end{equation*}

Then, by Theorem \ref{stability} and a diagonal argument (as $0 < T < + \infty$ is arbitrary) we find a global in time solution of the (MHDG) equations. \Endproof

\subsection{Solutions of the  advection-diffusion problem with initial data in  $L^2_{w_\gamma}$.}
Following essentially the same ideas of the proof of Theorem \ref{weightedMHD}, this result is  easily adapted for following advection-diffusion problem:

  \begin{Theorem}\label{advection} Within the hypothesis of Theorem  \ref{weightedMHD}, let $\vv, \vc$ be a time dependent divergence free vector-field such that, for every $T>0$, we have  $\vv, \vc \in L^3((0,T),L^3_{w_{3\gamma/2}})$. Then, the advection-diffusion problem
\begin{equation*}
(AD) \left\{ \begin{array}{ll}\vspace{2mm} 
\partial_t \vu = \Delta \vu - (\vv  \cdot \nabla) \vu + (\vc  \cdot \nabla) \vb - \nabla p + \vN \cdot \F,  \\ \vspace{2mm}
\partial_{t} \vb = \Delta \vb - (\vv  \cdot \nabla) \vb + (\vc  \cdot \nabla) \vu -\vN q + \vN \cdot \G 
, \\ \vspace{2mm}
\vN \cdot  \vu=0, \,  \vN \cdot \vb =0, \\ \vspace{2mm}
\vu(0,\cdot)=\vu_0, \, \vb(0,\cdot)=\vb_0, 
\end{array}
\right.     
\end{equation*} has a solution $(\vu, \vb, p , q )  $ which satisfies the statements of Theorem \ref{weightedMHD}.
 \end{Theorem}
\pv  For the initial data  $(\vu_{0,R}, \vb_{0,R})=(\P(\phi_R \vu_0), \P(\phi_R \vb_0))\in L^{2}(\mathbb{R}^3) $, for the forcing tensors $(\F_{R}, \G_R)= (\phi_R \F, \phi_R \G) \in L^2((0,T), L^2)$, and for $(\vv_R, \vc_R)=\mathbb{P}(\phi_R \vv, \phi_R \vc)$, with the arguments as in the proof of Proposition \ref{Prop-approx-leray}  we construct $(\vu_{R, \epsilon}, \vb_{R, \epsilon}, p_{R, \epsilon}, q_{R, \epsilon})$, a solution  of the approximated system 
\begin{equation*}
\left\{ \begin{array}{ll}\vspace{2mm} 
\partial_t \vu_{R, \epsilon} = \Delta \vu_{R, \epsilon} - ((\vv_{R} * \theta_\epsilon ) \cdot \nabla) \vu_{R, \epsilon} + ((\vc_{R} * \epsilon) \cdot \nabla) \vb_{R, \epsilon} - \nabla p_{R, \epsilon} + \vN \cdot \F_R,  \\ \vspace{2mm}
\partial_{t} \vb_{R, \epsilon} = \Delta \vb_{R, \epsilon} - ((\vv_{R} * \epsilon) \cdot \nabla) \vb_{R, \epsilon} + ((\vc_{R} * \epsilon) \cdot \nabla) \vu_{R, \epsilon} -\nabla q_{R, \epsilon} + \vN  \cdot \G_R 
, \\ \vspace{2mm}
\vN \cdot  \vu_{R, \epsilon}=0, \,  \vN \cdot \vb_{R, \epsilon} =0, \\ \vspace{2mm}
\vu_{R, \epsilon}(0,\cdot)=\vu_{0,R}, \, \vb_{R, \epsilon}(0,\cdot)=\vb_{0,R},  
\end{array}
\right.     
\end{equation*}
where $ (\vu_{R, \epsilon}, \vb_{R, \epsilon}) \in \mathcal{C} ([0, T),L^{2}(\Rt))\cap L^{2}([0, T), \dot{H}^{1}(\Rt))$ for every $0<T<+\infty$, and $(\vu_{R, \epsilon}, \vb_{R, \epsilon}, p_{R, \epsilon}, q_{R, \epsilon})$ verify all  the assumptions of Theorem  \ref{estimates} (with energy equality).
We remark that 
\begin{align}
\label{cvc}
\begin{split}
    \| (\vv_{R}*\theta_\epsilon,\vc_{R}*\theta_\epsilon)\|_{L^3((0,T), L^3_{w_{3\gamma/2}} ) } 
    &\leq \| (\mathcal{M}_{\vv_R}, \mathcal{M}_{\vc_R} )\|_{L^3((0,T), L^3_{w_{3\gamma/2}} ) } \\
    &\leq C_\gamma       \| (\vv, \vc )\|_{L^3((0,T), L^3_{w_{3\gamma/2}} ) }.
\end{split}
\end{align}
 
We write  $(\vu_{0,n},\vb_{0,n})=(\vu_{0,R_n}, \vb_{0,R_n} )$, $ (  \mathbb{F}_n, \G_n ) = (\mathbb{F}_{R_n}, \mathbb{G}_{R_n} )$, $(\vv_n,\vc_n) =(\vv_{R_n} * \epsilon_n, \vv_{R_n} * \epsilon_n)$ and $(\vu_n, \vb_n)=(\vu_{R_n,\epsilon_n},\vb_{R_n,\epsilon_n})$.  

As  $(\vu_{0,n},\vb_{0,n})$ is strongly convergent to $(\vu_{0}, \vb_0)$ in $L^2_{w_\gamma}$,   $(\mathbb{F}_n, \G_n)$ is strongly convergent to $(\mathbb{F}, \G)$  in $L^2((0,T), L^2_{w_\gamma})$,  and  moreover,  as $(\vv_n, \vc_n)$ is bounded in $L^3((0,T),   L^3_{w_{3\gamma/2}})$ (since \eqref{cvc}), we can apply Theorem \ref{stability}, then  there exist $(\vu, \vb, {\bf V},  {\bf C}, p,q)$ and there exists an increasing  sequence $(n_k)_{k\in\mathbb{N}}$ with values in $\mathbb{N}$ such that:

  \begin{itemize} 
 \item[$\bullet$] $(\vu_{n_k}, \vb_{n_k} )$ converges *-weakly to $(\vu, \vb)$ in $L^\infty((0,T_0), L^2_{w_\gamma})$, $ (\vN\vu_{n_k}, \vN\vb_{n_k} )$ converges weakly to $\vN (\vu, \vb)$ in $L^2((0,T_0),L^2_{w_\gamma})$.
 \item[$\bullet$] $(\vv_{n_k}, \vc_{n_k} )$ converges weakly to $({\bf V}, {\bf C})$ in $L^3((0,T_0), L^3_{w_{3\gamma/2}})$, the sequence $p_{n_k}$ converges weakly to $p$ in $L^{3}((0,T_0),L^{6/5}_{w_{\frac {6\gamma}5}})+L^{2}((0,T_0),L^{2}_{w_\gamma})$ and similarly for $q_{n_k}$.
  \item[$\bullet$] $(\vu_{n_k}, \vb_{n_k} )$ converges strongly  to $(\vu, \vb)$ in  $L^2_{\rm loc}([0,T_0)\times\mathbb{R}^3)$, \end{itemize}
 and moreover,  $(\vu, \vb, {\bf V},  {\bf C}, p,q)$ is  a solution of the   advection-diffusion problem 
 \begin{equation*}
 \left\{ \begin{array}{ll}\vspace{2mm} 
\partial_t \vu = \Delta \vu - ({\bf V} \cdot \nabla) \vu + (\vc \cdot \nabla) \vb - \nabla p + \vN \cdot \F,  \\ \vspace{2mm}
\partial_{t} \vb = \Delta \vb - ({\bf V} \cdot \nabla) \vb + ({\bf C} \cdot \nabla) \vu -\vN q + 
\vN  \cdot \G 
, \\ \vspace{2mm}
\vN \cdot  \vu=0, \,  \vN \cdot \vb =0, \\ \vspace{2mm}
\vu(0,\cdot)=\vu_0, \, \vb(0,\cdot)=\vb_0.  
\end{array}
\right.     
\end{equation*}
which verifies:
 \begin{itemize}  
 \item[$\bullet$] the map $t\in [0,T_0)\mapsto (\vu(t), \vb(t))$ is weakly continuous from $[0,T_0)$ to $L^2_{w_\gamma}$, and is strongly continuous at $t=0$.
 \item[$\bullet$]  there exists a non-negative locally finite measure $\mu$ on $(0,T)\times\mathbb{R}^3$ such that we have the local energy balance (\ref{energloc1}).   
 \end{itemize} 
If we verify that   ${\bf V}=\vv$ and ${\bf C}=\vc$ the proof is finished. As $\vv_n=\theta_{\epsilon_n}*(\vv_n-\vv)+ \theta_{\epsilon_n}*\vv$, and $\vc_n=\theta_{\epsilon_n}*(\vc_n-\vc)+ \theta_{\epsilon_n}*\vc$  then we get that  $(\vv_{n_k},\vc_{n_k})$ are strongly   convergent to $(\vv, \vc)$ in   $L^3_{\rm loc}([0,T_0)\times\mathbb{R}^3)$, hence we have  ${\bf V}=\vv$ and ${\bf C}=\vc$.
\Endproof

\section{Discretely self-similar suitable solutions for 3D MHD equations}\label{sec:self-similar}
In this section we give a proof of Theorem \ref{selfsimilar}. We fix $ 1 < \lambda < +\infty.  $
 \subsection{The linear problem.}
Let $\theta$ be a  non-negative and radially decreasing  function in $\mathcal{D}(\mathbb{R}^3)$ with $\int\theta\, dx=1$; We define $\theta_{\epsilon,t}(x)=\frac 1 {(\epsilon\sqrt t)^3}\ \theta(\frac x{\epsilon\sqrt t})$.  In order to study the mollified problem
\begin{equation*}
(MHD_\epsilon) \left\{ \begin{array}{ll}\vspace{2mm} 
\partial_t \vu_\epsilon = \Delta \vu_\epsilon - ((\vu_\epsilon * \theta_{\epsilon, t}) \cdot \nabla) \vu_\epsilon + ((\vb_\epsilon * \theta_{\epsilon, t}) \cdot \nabla) \vb_\epsilon - \nabla p + \vN \cdot \F,  \\ \vspace{2mm}
\partial_{t} \vb_\epsilon = \Delta \vb_\epsilon - ((\vu_\epsilon  *\theta_{\epsilon, t}) \cdot \nabla) \vb_\epsilon + ((\vb_\epsilon *\theta_{\epsilon, t}) \cdot \nabla) \vu_\epsilon -\vN q + \vN \cdot \G 
, \\ \vspace{2mm}
\vN \cdot  \vu_\epsilon =0, \,  \vN \cdot \vb_\epsilon =0, \\ \vspace{2mm}
\vu(0,\cdot)=\vu_0, \, \vb(0,\cdot)=\vb_0,  
\end{array}
\right.     
\end{equation*}
we consider the linearized problem
\begin{equation*}
(LMHD) \left\{ \begin{array}{ll}\vspace{2mm} 
\partial_t \vu = \Delta \vu - ((\vv * \theta_{\epsilon, t}) \cdot \nabla) \vu + ((\vc * \theta_{\epsilon, t}) \cdot \nabla) \vb - \nabla p + \vN \cdot \F,  \\ \vspace{2mm}
\partial_{t} \vb = \Delta \vb - ((\vv  *\theta_{\epsilon, t}) \cdot \nabla) \vb + ((\vc *\theta_{\epsilon, t}) \cdot \nabla) \vu -\vN q + \vN \cdot \G 
, \\ \vspace{2mm}
\vN \cdot  \vu=0, \,  \vN \cdot \vb =0, \\ \vspace{2mm}
\vu(0,\cdot)=\vu_0, \, \vb(0,\cdot)=\vb_0.  
\end{array}
\right.  
\end{equation*}

  \begin{Lemma}\label{solvLNS}
Let $1<\gamma\leq 2$. Let $\vu_{0}, \vb_{0}$ be a $\lambda$-DSS divergence-free vector fields which belong to  $L^2_{w_\gamma}(\mathbb{R}^3)$. Let $\mathbb{F}, \G$ be a $\lambda$-DSS tensors  wich satisfies $\mathbb{F}, \G \in L^2_{loc} ((0,+\infty), L^2_{w_\gamma})$. Moreover, let $\vv, \vc$ be a $\lambda$-DSS  time-dependent divergence free vector-field such that for every $T>0$,  $\vv, \vc \in L^3((0,T),L^3_{w_{3\gamma/2}})$. \\

Then, the linearized advection-diffusion problem (LMHD) has a unique solution $(\vu, \vb, p , q )$ which satisfies all the conclusions of Theorem \ref{advection}. Moreover, the functions $\vu, \vb  $ are $\lambda$-DSS vector fields.
 \end{Lemma}
  
\pv  As we have $\vert \vv(t,.)*\theta_{\epsilon,t}\vert\leq \mathcal{M}_{\vu(t,.)} $ then we can write
 $$ \|( \vv(t)*\theta_{\epsilon,t}, \vc(t)*\theta_{\epsilon,t} )\|_{L^3((0,T), L^3_{w_{3\gamma/2}})}\leq C_\gamma \|( \vv, \vc )\|_{L^3((0,T), L^3_{w_{3\gamma/2}})}.$$ 
 
Theorem \ref{advection} gives solution $(\vu, \vb, p , q )$ in the interval of time $(0,T)$. Moreover, as $\vu*\theta_{\epsilon,t}, \vb*\theta_{\epsilon,t}$  belong the space  to $L^2_t L^\infty_x(K)$ for every compact subset $K$ of  $(0,T)\times\mathbb{R}^3 $, we can use Corollary \ref{unique} to conclude that this solution $(\vu, \vb, p , q )$ is unique.\\

 We will prove that this solution is $\lambda$-DSS. Let $\tilde \vu (t,x)=\frac 1 \lambda \vu(\frac t{\lambda^2}, \frac x \lambda)$ and $\tilde \vb (t,x)=\frac 1 \lambda \vb(\frac t{\lambda^2}, \frac x \lambda)$. Remark that $(\vv*\theta_{\epsilon,t}$ and $ \vc*\theta_{\epsilon,t})$ are  $\lambda$-DSS and then we get  $(\tilde \vu, \tilde \vb, \tilde p, \tilde q)$, where $\tilde p$ and$ \tilde q$ are always defined through the obvious formula, is a solution of $(LMHD_\epsilon)$ on  $(0,T)$. Thus, we have the identities $(\tilde \vu, \tilde \vb, \tilde p, \tilde q) = (\vu, \vb, p , q )$ from which we conclude that $(\vu, \vb, p , q )$ are $\lambda$-DSS. 
 
 \Endproof

 \subsection{The mollified Navier--Stokes equations.}
 
 For $\vv, \vc \in L^3((0,T), L^3_{w_{3\gamma/2}})$ the terms
 $\vu, \vb$ of the solution provided by Lemma \ref{solvLNS} belongs to $L^3((0,T), L^3_{w_{3\gamma/2}})$ by interpolation. Then the map $L_\epsilon : (\vv, \vc)  \mapsto (\vu, \vb) $ where $L_\epsilon(\vv, \vc)=(\vu, \vb)$ is well defined from
 $$ X_{T,\gamma}=\{ (\vv, \vc) \in L^3((0,T), L^3_{w_{3\gamma/2}})\ /\ \vb \text{ is } \lambda-\text{DSS}\}$$ to $X_{T,\gamma}$. 
At this point, we introduce the following technical lemmas: 
\begin{Lemma} For $4/3<\gamma$, $X_{T,\gamma}$ is a Banach space for the equivalent norms $\|(\vv,\vc) \|_{L^3((0,T),L^3_{w_{3\gamma/2}})}$ and $\|(\vv,\vc)\|_{L^3((0,T/{\lambda^2}),\times B(0,\frac 1 \lambda))}$. \end{Lemma}
For a proof of this result see the  Lemma 12 in \cite{PF_PG}.
 
 \begin{Lemma}\label{compact} For $4/3<\gamma\leq 2$, the mapping $L_\epsilon$ is continuous and compact on $X_{T,\gamma}$.
 \end{Lemma} 
 
 \pv Let $(\vv_n,\vc_n)$ be a bounded sequence in $X_{T,\gamma}$ and let $(\vu_n,\vb_n)=L_\epsilon(\vv_n,\vc_n)$. Remark that the sequence $(\vv_n(t)*\theta_{\epsilon,t},\vc_n(t)*\theta_{\epsilon,t} )$ is bounded in $X_{T,\gamma}$ and then by Theorem \ref{estimates} and Corollary \ref{passive} we have that the sequence $(\vu_n,\vb_n)$ is bounded in $L^\infty((0,T), L^2_{w_\gamma})$ and  moreover $(\vN\vu_n, \vN\vb_n)$ is bounded in $L^2((0,T),L^2_{w_\gamma})$.\\
 
Thus, by Theorem \ref{stability}  there exists $\vu_\infty$, $\vb_\infty$, $p_\infty$, $q_\infty$, ${\bf  V}_\infty$, ${\bf  C}_\infty$   and an increasing  sequence $(n_k)_{k\in\mathbb{N}}$ with values in $\mathbb{N}$ such that we have:
  \begin{itemize} 
 \item[$\bullet$] $(\vu_{n_k}, \vb_{n_k})$ converges *-weakly to $(\vu_\infty, \vb_\infty)$ in $L^\infty((0,T), L^2_{w_\gamma})$, $ (\vN\vu_{n_k}, \vN\vb_{n_k})$ converges weakly to $(\vN\vu_\infty, \vN\vb_\infty) $ in $L^2((0,T),L^2_{w_\gamma})$.
 \item[$\bullet$] $( \vv_{n_k}*\theta_{\epsilon,t}, \vc_{n_k}*\theta_{\epsilon,t} )$ converges weakly to $({\bf V}_\infty,  {\bf C}_\infty)$ in $L^3((0,T), L^3_{w_{3\gamma/2}})$. 
 \item[$\bullet$]  The  terms $(p_{n_k}, q_{n_k})$ converge weakly to $(p_\infty, q_\infty)$ in $L^{3}((0,T),L^{6/5}_{w_{\frac {6\gamma}5}})+L^{2}((0,T),L^{2}_{w_\gamma})$.
  \item[$\bullet$] $(\vu_{n_k}, \vb_{n_k})$ converges strongly  to $(\vu_\infty, \vb_\infty)$ in  $L^2_{\rm loc}([0,T)\times\mathbb{R}^3)$ : for every  $T_0\in (0,T)$ and every $R>0$, we have
  $$\lim_{k\rightarrow +\infty} \int_0^{T_0} \int_{\vert y\vert<R} \vert \vu_{n_k}(s,y)-\vu_\infty(s,y)\vert^2 + \vert \vb_{n_k}(s,y)-\vb_\infty(s,y)\vert^2 \, ds\, dy=0,$$
  \item[$\bullet$] and
    \begin{equation*}
    \left\{ \begin{array}{ll}\vspace{2mm} 
    \partial_t \vu_\infty = \Delta \vu_\infty - (\vv_\infty \cdot \nabla) \vu_\infty + (\vc_\infty \cdot \nabla) \vb_\infty - \nabla {p}_\infty + \vN \cdot \F,  \\ \vspace{2mm}
    \partial_{t} \vb_\infty = \Delta \vb_\infty - (\vv_\infty \cdot \nabla) \vb_\infty + (\vc_\infty \cdot \nabla) \vu_\infty -\vN {q}_\infty + 
    \vN  \cdot \G 
    , \\ \vspace{2mm}
    \vN \cdot  \vu_\infty=0, \,  \vN \cdot \vb_\infty =0,\\ \vspace{2mm}
     \vu_{0,\infty}=\vu_0, \,  \vb_{0,\infty} =\vb_0.
    \end{array}
    \right.     
    \end{equation*}
  \end{itemize}
 We will prove the compactness of  $L_\epsilon$. As before $\sqrt{w_\gamma}\vv_n$  is bounded in $L^{10/3}((0,T)\times\mathbb{R}^3)$ by interpolation hence strong convergence of $(\vu_{n_k},\vb_{n_k})$ in $L^2_{\rm loc}([0,T)\times\mathbb{R}^3)$ implies the strong convergence of $(\vu_{n_k},\vb_{n_k})$ in $L^3_{\rm loc}((0,T)\times\mathbb{R}^3)$.
  
  Moreover, we have that  $(\vu_\infty,\vb_\infty)$ is still $\lambda$-DSS (a property that is stable under weak limits). With these information we obtain that $\vu_\infty,\vb_\infty \in X_{T,\gamma}$ and we have
  $$ \lim_{n_k\rightarrow +\infty} \int_0^{\frac T {\lambda^2}} \int_{B(0,\frac 1 \lambda)} \vert \vv_{n_k}(s,y)-\vv_\infty(s,y)\vert^3\, ds\, dy=0,$$ which  proves that $L_\epsilon$ is compact.\\
  
 To finish this proof, we prove the continuity of $L_\epsilon$. Let $(\vv_n, \vc_n)$ be such that  $(\vv_n, \vc_n)$ is convergent to $(\vv_\infty, \vc_\infty)$ in $X_{T,\gamma}$. Then  we have ${\bf V}_\infty=\vv_\infty*\theta_{\epsilon,t}$, ${\bf C}_\infty=\vc_\infty*\theta_{\epsilon,t}$, and $\vu_\infty= L_\epsilon(\vv_\infty, \vc_\infty)$, and thus,  the relatively compact sequence $(\vu_n,\vb_n)$ can have only one limit point. In conclusion, it must be convergent and this  proves that $L_\epsilon$ is continuous.\Endproof

 \begin{Lemma}\label{apriori} Let $4/3<\gamma\leq 2$. If $\mu\in [0,1]$ and  $(\vu, \vb)$ solves $(\vu, \vb)=\mu L_\epsilon(\vu, \vb)$ then $$\|(\vu, \vb)\|_{X_{T, \gamma}}\leq C_{\vu_0,\mathbb{F},\gamma, T, \lambda}$$ where the constant $C_{\vu_0,\mathbb{F},\gamma, T, \lambda}$ depends only on $\vu_0$, $\mathbb{F}$, $\gamma$, $T$ and $\lambda$ (but not on $\mu$ nor on $\epsilon$).
 \end{Lemma}
 
 \pv We let $(\vu, \vb)= (\mu \tilde \vu,\mu \tilde \vb)$, so that 
\begin{equation*}
 \left\{ \begin{array}{ll}\vspace{2mm} 
\partial_t \tilde \vu = \Delta \tilde \vu - ((\vu * \theta_{\epsilon, t}) \cdot \nabla) \tilde \vu + ((\vb * \theta_{\epsilon, t}) \cdot \nabla) \tilde \vb - \nabla p + \vN \cdot \F,  \\ \vspace{2mm}
\partial_{t} \tilde \vb = \Delta \tilde \vb - ((\vu  *\theta_{\epsilon, t}) \cdot \nabla) \tilde \vb + ((\vb *\theta_{\epsilon, t}) \cdot \nabla) \tilde \vu -\vN q + \vN \cdot \G 
, \\ \vspace{2mm}
\vN \cdot  \tilde \vu=0, \,  \vN \cdot \tilde \vb =0, \\ \vspace{2mm}
\tilde \vu(0,\cdot)=\vu_0, \, \tilde \vb(0,\cdot)=\vb_0.  
\end{array}
\right.     
\end{equation*}
 
 Multiplying by $\mu$, we find that
\begin{equation*}
 \left\{ \begin{array}{ll}\vspace{2mm} 
\partial_t \vu = \Delta \vu - ((\vu * \theta_{\epsilon, t}) \cdot \nabla) \vu + ((\vb * \theta_{\epsilon, t}) \cdot \nabla) \vb - \nabla (\mu p) + \vN \cdot \mu \F,  \\ \vspace{2mm}
\partial_{t} \vb = \Delta \vb - ((\vu  *\theta_{\epsilon, t}) \cdot \nabla) \vb + ((\vb *\theta_{\epsilon, t}) \cdot \nabla) \vu -\vN (\mu q) + \vN \cdot \mu \G 
, \\ \vspace{2mm}
\vN \cdot  \vu=0, \,  \vN \cdot \vb =0, \\ \vspace{2mm}
\vu(0,\cdot)=\mu \vu_0, \, \vb(0,\cdot)=\mu \vb_0.  
\end{array}
\right.     
\end{equation*}
 
 Corollary \ref{active} allows us to take $T_0\in (0,T)$ such that
 $$ C_\gamma  \left(1+\|(\vu_0,\vb_0)\|_{L^2_{w_\gamma}}^2+\int_0^{T_0} \|(\mathbb{F}, \G)\|_{L^2_{w_\gamma}}^2\, ds\right)^2\, T_0\leq 1,$$ which implies
 $$ C_\gamma  \left(1+\| \mu (\vu_0,\vb_0) \|_{L^2_{w_\gamma}}^2+\int_0^{T_0} \ \| \mu (\mathbb{F}, \G) \|_{L^2_{w_\gamma}}^2\, ds\right)^2\, T_0\leq 1.$$
 Then we have the controls
 \begin{align*}
     \begin{split}
    \sup_{0\leq t\leq T_0} & \| (\vu, \vb)(t)\|_{L^2_{w_\gamma}}^2 + { \int_0^{T_0} \|\vN(\vu, \vb)\|_{L^2_{w_\gamma}}^2\, ds } \\
         &\leq C_\gamma (1 + \mu^2  \|(\vu_0,\vb_0)\|_{L^2_{w_\gamma}}^2 +\mu^2\int_0^{T_0} \|(\mathbb{F}, \G)\|_{L^2_{w_\gamma}}^2\, ds ).
     \end{split}
 \end{align*}
 In particular, by interpolation
 $$ \int_0^{T_0} \|(\vu, \vb)\|_{L^3_{w_{3\gamma/2}}}^3\, ds
 $$ 
 is bounded by a constant $C_{\vu_0,\mathbb{F},\gamma, T}$ and we can go back from $T_0$ to $T$, using the self-similarity property.
 \Endproof

 \begin{Lemma}\label{schauder} Let $4/3<\gamma\leq 2$. There is at least one solution $(\vu_\epsilon, \vb_\epsilon)$ of the problem $(\vu_\epsilon, \vb_\epsilon)=L_\epsilon(\vu_\epsilon, \vb_\epsilon)$.
 \end{Lemma}

\pv  The uniform a priori estimates for the fixed points of $\mu L_\epsilon$ for $0\leq\mu\leq 1$ given by Lemma \ref{apriori} and Lemma \ref{compact} permit to apply Leray--Schauder principle and Schaefer theorem.\Endproof 

\subsection{Proof of Theorem \ref{selfsimilar}.}
We consider $(\vu_\epsilon, \vb_\epsilon)$ solutions of $(\vu_\epsilon, \vb_\epsilon)=L_\epsilon(\vu_\epsilon, \vb_\epsilon)$ given by Lemma \ref{schauder}.

By Lemma \ref{apriori} and Lemma \ref{indt}, we have  $\vu_\epsilon*\theta_{\epsilon,t},\vb_\epsilon*\theta_{\epsilon,t}$ are bounded  in $L^3((0,T), L^3_{w_{3\gamma/2}})$. Theorem \ref{estimates} and Corollary \ref{passive} allows us to conclude that $\vu_\epsilon, \vb_\epsilon$ are bounded in $L^\infty((0,T), L^2_{w_\gamma})$ and $\vN \vu_\epsilon, \vN \vb_\epsilon$ are bounded in $L^2((0,T),L^2_{w_\gamma})$.
 
Theorem \ref{stability} gives $\vu$,  $\vb$, $p$, $q$,  ${\bf  v}$ and ${\bf c}$   and a decreasing  sequence $(\epsilon_k)_{k\in\mathbb{N}}$ converging to $0$, such that
  \begin{itemize} 
 \item[$\bullet$] $ (\vu_{\epsilon_k},\vb_{\epsilon_k} ) $ converges *-weakly to $(\vu, \vb)$ in $L^\infty((0,T), L^2_{w_\gamma})$, $(\vN\vu_{\epsilon_k},\vN\vb_{\epsilon_k} )$ converges weakly to $( \vN\vu, \vN\vu )$ in $L^2((0,T),L^2_{w_\gamma})$
 \item[$\bullet$] $( \vu_{\epsilon_k}*\theta_{\epsilon_k,t}, \vb_{\epsilon_k}*\theta_{\epsilon_k,t} )$ converges weakly to $(\vv, \vc)$ in $L^3((0,T), L^3_{w_{3\gamma/2}})$
  \item[$\bullet$]  the associated pressures $p_{\epsilon_k}$ and $ q_{\epsilon_k} $ converge weakly to $p$ and $q$ in $L^{3}((0,T),L^{6/5}_{w_{\frac {6\gamma}5}})+L^{2}((0,T),L^{2}_{w_\gamma})$
  \item[$\bullet$] $(\vu_{\epsilon_k},\vb_{\epsilon_k} )$ converges strongly  to $( \vu, \vb)$ in  $L^2_{\rm loc}([0,T)\times\mathbb{R}^3)$
  \item[$\bullet$] and
    \begin{equation*}
    \left\{ \begin{array}{ll}\vspace{2mm} 
    \partial_t \vu = \Delta \vu - (\vv \cdot \nabla) \vu + (\vc \cdot \nabla) \vb - \nabla {p} + \vN \cdot \F,  \\ \vspace{2mm}
    \partial_{t} \vb = \Delta \vb - (\vv \cdot \nabla) \vb + (\vc \cdot \nabla) \vu -\vN {q} + 
    \vN  \cdot \G 
    , \\ \vspace{2mm}
    \vN \cdot  \vu=0, \,  \vN \cdot \vb =0,\\  \vspace{2mm}
     \vu_{0}=\vu_0, \,  \vb_{0} =\vb_0,
    \end{array}
    \right.     
    \end{equation*}
  \end{itemize}
  The proof is finished if $\vv=\vu$ and $\vc=\vb$.  As we have  $\vu_{\epsilon_k}*\theta_{\epsilon_k,t}=(\vu_{\epsilon_k}-\vu)*\theta_{\epsilon_k,t}+ \vu *\theta_{\epsilon_k,t}$. We just need to remark that $\vu*\theta_{\epsilon,t}$ converges strongly in $L^2_{\rm loc}((0,T)\times \mathbb{R}^3)$ as $\epsilon$ goes to $0$ (we use dominated convergence as it is bounded by $\mathcal{M}_\vu$ and converges strongly to $\vu$ in $L^2_{\rm loc}(\mathbb{R}^3)$ for each fixed $t$ ) and $ \vert (\vu-\vu_\epsilon)*\theta_{\epsilon,t}\vert\leq \mathcal{M}_{\vu-\vu_\epsilon}$. In a similar way we prove $\vc=\vb$. \Endproof
  
 \appendix

\section{Approximated system}
Let $\theta \in \mathcal{D}(\Rt)$ be a non-negative, radial and radially decreasing function such that $\int_{\Rt}\theta(x)dx=1$. For $\varepsilon >0$ we let $\theta_\varepsilon (x) = \frac{1}{\varepsilon ^3} \theta (\frac{x}{\varepsilon})$.
\begin{Proposition}\label{Prop-approx-leray}  Let $\vu_0 \in L^2(\Rt)$, $\vb_0 \in L^{2}(\Rt)$ be divergence free vector fields. Let $\mathbb{F } = (F_{i,j})_{1 \leq i,j \leq 2 }  $ and $\mathbb{G } = (G_{i,j})_{1 \leq i,j \leq 2 }$ be tensor forces such that $\mathbb{F }
, \mathbb{G }  
\in L^{2}((0, T), L^2)$, for all $T< T_\infty$. 

Then there exists a unique solution $(\vue, \vbe, p_\varepsilon , q_\varepsilon)$ of the following approximated system
\begin{equation*}
(MHDG_\varepsilon ) \left\{ \begin{array}{ll}\vspace{2mm} 
\partial_t \vu = \Delta \vu - [(\vu\ast \theta_{\varepsilon}) \cdot \nabla] \vu + [(\vb \ast \theta_{\varepsilon}) \cdot \nabla] \vb - \nabla p + \vN \cdot \F,  \\ \vspace{2mm}
\partial_{t} \vb = \Delta \vb - [(\vue\ast \theta_{\varepsilon}) \cdot \nabla] \vb + [(\vb \ast \theta_\varepsilon) \cdot \nabla] \vu - \vN q + 
\vN  \cdot \G 
,\\ \vspace{2mm}
\vN \cdot  \vu=0, \,  \vN \cdot \vb =0, \\ \vspace{2mm}
\vu(0,\cdot)=\vu_0, \, \vb(0,\cdot)=\vb_0, 
\end{array}
\right.     
\end{equation*}
on $[0,T_\infty)$ such that:
\begin{enumerate}
\item[$\bullet$] $\vu_\varepsilon, \vb_\varepsilon \in L^{\infty} ([0, T),L^{2}(\Rt))\cap L^{2}([0, T), \dot{H}^{1}(\Rt))$, $ p_\varepsilon ,q_\varepsilon \in L^2((0,T), \dot H ^{-1}) +L^2((0,T), L^2) $, for all $0 <T < T_\infty $

\item[$\bullet$]  the pressure $p_\varepsilon$ and the term $q_\varepsilon$ are related to $\vu_\varepsilon, \vb_\varepsilon$, $\mathbb{F}$ and $\mathbb{G}$ by 
$$\pe = \sum_{1 \leq i,j\leq 3} \mathcal{R}_i \mathcal{R}_j ((u_{\varepsilon, i}*\te) u_{\varepsilon,j} - (b_{\varepsilon,i}*\te) b _{\varepsilon,j} - F_{i,j} ),$$
and
$$\qe =  \sum_{1 \leq i,j\leq 3} \mathcal{R}_i \mathcal{R}_j ( [(u_{\varepsilon,i}*\te) b_{\varepsilon,j} - (b_{\varepsilon,j}*\te) u _{\varepsilon,i}]- G_{ij}),$$ where $\Ri_i =\frac{\partial_i}{\sqrt{-\Delta}}$ denote always the Riesz transforms. In particular, $ p_\varepsilon ,q_\varepsilon \in  L^{4}((0,T),L^{6/5}) + L^{2}((0,T),L^{2})$.
\item[$\bullet$]The functions $(\vu_\varepsilon, \vb_\varepsilon, \F, \G)$ verify the following global energy equality: 
\begin{equation}\label{eneqpo}\begin{split}
    \partial_t(\frac {\vert\vu_\varepsilon \vert^2 
    + |\vbe|^2 }2) =&\Delta(\frac {\vert\vu_\varepsilon \vert^2 + |\vbe|^2 }2)-\vert\vN \vu_\varepsilon \vert^2 - |\vN \vbe|^2 \\
    &- \vN\cdot\left( (\frac{\vert\vu_\varepsilon \vert^2}2 + \frac{\vert\vb_\varepsilon \vert^2}2)(\vu_\varepsilon * \theta_\varepsilon ) + \pe \vue  \right) \\
    & + \vN \cdot ((\vu_\varepsilon \cdot \vb_\varepsilon) (\vb_\varepsilon * \theta_\varepsilon) + \qe \vbe) \\
    & + \vu_\varepsilon \cdot(\vN\cdot\mathbb{F})+ \vbe  \cdot(\vN\cdot\mathbb{G}). 
\end{split}
\end{equation}
and
\begin{eqnarray*}
\begin{split}
    \Vert \vu_\varepsilon (t) \Vert^{2}_{L^2} & + \Vert \vb_\varepsilon (t) \Vert^{2}_{L^2} + 2 \int_{a}^{t} (\Vert \vN \vu_\varepsilon (s) \Vert^{2}_{L^2} + \Vert \vN \vb_\varepsilon (s) \Vert^{2}_{L^2} )	 ds \nonumber\\
    =&  \Vert \vu_\varepsilon (a) \Vert^{2}_{L^2} + \Vert \vb_\varepsilon  (a) \Vert^{2}_{L^2} 
     \nonumber \\
     & + \sum_{1 \leq i,j \leq 3} (  \int_{a}^{t}  \int   F_{i,j}\partial_iu_{_\varepsilon,j}\    \, dx\, ds + \int_{a}^{t}  \int   G_{i,j}\partial_ib_{_\varepsilon,j}\    \, dx\, ds ) ,
\end{split}
\end{eqnarray*}
which implies in particular
\begin{equation}
\label{eneq1}
\begin{split}
    \Vert \vu_\varepsilon(t) \Vert^{2}_{L^2} &+ \Vert \vb_\varepsilon(t) \Vert^{2}_{L^2} + \int_{0}^{t} (\Vert \vN \vu_\varepsilon(s) \Vert^{2}_{L^2} + \Vert \vN \vb_\varepsilon(s) \Vert^{2}_{L^2} )	 ds \nonumber\\
    & \leq  \Vert \vu_0 \Vert^{2}_{L^2} +  \Vert \vb_0 \Vert^{2}_{L^2} 
    +   c (\Vert \F \Vert^{2}_{L^{2}_{t} L^{2}_{x}} + \Vert \G \Vert^{2}_{L^{2}_{t} L^{2}_{x}}).   
\end{split}
\end{equation}
\end{enumerate}
\end{Proposition}

\pv We consider $0 < T < T_1 < T_\infty$ and the  space $E_T= \mathcal{C}([0,T], L^2(\Rt))  \cap L^2((0, T) \dot{H}^{1}(\Rt))$ doted with the norm $\Vert \cdot \Vert_{T}= \Vert \cdot \Vert_{L^{\infty}_{t}L^{2}_{x}}+ \Vert \cdot \Vert_{L^{2}_{t}\dot{H}^{1}_{x}}$. We will construct simultaneously $\vue$ and $\vbe$. For this we will consider the space $E_T \times E_T$ with the norm $\Vert (\vue, \vbe) \Vert_{T}= \Vert \vue \Vert_{T}+ \Vert \vbe \Vert_{T}$.

We use the Leray projection operator in order to express the problem $(MHDG_\varepsilon)$ in terms of a fixed point problem. We let 
\begin{equation*}
    a= e^{ t \Delta  } (\vv_0, \vc_0)+ \int_{0}^{t} e^{(t-s)\Delta}  \P( \vN \cdot \F,  \vN \cdot \G)(s,\cdot) ds
\end{equation*}
and
\begin{equation*}
    B((\vu,\vb), (\vv,\vc)) = ( \, B_1((\vu,\vb)\, , \,  (\vv,\vc)), B_2((\vu,\vb), (\vv,\vc)) \, ) ,
\end{equation*}
where
\begin{eqnarray*}\label{eq01} \nonumber
B_1((\vu,\vb), (\vv,\vc)) = & &
 \int_{0}^{t} e^{(t-s)\Delta }  \P ([(\vu \ast \theta_{\varepsilon}) \cdot \nabla) \vv -  [(\vv \ast \theta_{\varepsilon}] \cdot \nabla] \vc )(s,\cdot) ds ,
\end{eqnarray*}
\begin{eqnarray*}\label{eq01} \nonumber
B_2((\vu,\vb), (\vv,\vc)) = & &
 \int_{0}^{t} e^{(t-s)\Delta }  \P ([(\vu \ast \theta_{\varepsilon}) \cdot \nabla] \vc -  [(\vb \ast \theta_{\varepsilon}) \cdot \nabla] \vv )(s,\cdot) ds .
\end{eqnarray*}
Then $$(\vue, \vbe, \pe, \qe) \in E_T^2 \times \left( L^2((0,T), \dot H ^{-1}) +L^2((0,T), L^2) \right)^2 $$ is a solution of $(MHDG_\varepsilon )$ if and only if $(\vue, \vbe)$ is a fixed point for the application $(\vu,\vb) \mapsto a + B((\vu,\vb),(\vu,\vb))$ and 
$$\pe  = \sum_{1 \leq i,j\leq 3} \mathcal{R}_i \mathcal{R}_j ((u_{\varepsilon, i}*\te) u_{\varepsilon,j} - (b_{\varepsilon,i}*\te) b _{\varepsilon,j} - F_{i,j} ),$$
and
$$\qe  = \sum_{1 \leq i,j\leq 3} \mathcal{R}_i \mathcal{R}_j ( [(u_{\varepsilon,i}*\te) b_{\varepsilon,j} - (b_{\varepsilon,j}*\te) u _{\varepsilon,i}]- G_{ij}).$$

We will use the Piccard's point fixed theorem. In order to study  the linear terms, recall the following estimates, for a proof see   \cite{LR16}, Theorem $12.2$, page $352$.
\begin{Lemma} Let $f \in L^{2}(\Rt)$ and $g \in L^{2}_{t}\dot{H}^{-1}_{x}$. We have: 
\begin{enumerate}
\item[1)] $\Vert e^{ t \Delta } f \Vert_{T} \leq c \Vert f \Vert_{L^2}$. 
\item[2)] $\left\Vert \int_{0}^{t} e^{(t-s)\Delta} g(s,\cdot)ds \right\Vert_{T} \leq c(1+\sqrt{T}) \Vert g \Vert_{L^{2}_{t} \dot{H}^{-1}_{x}}$. 
\end{enumerate}		
\end{Lemma}	
By this lemma we have 
\begin{equation}\label{eq02}
\Vert  e^{ t \Delta  } (\vu_0, \vb_0) \Vert_{T} \leq c (\Vert \vu_0 \Vert_{L^2} + \Vert \vb_0 \Vert_{L^2}),
\end{equation}
and 
\begin{eqnarray}\label{eq03} \nonumber 
& & \left\Vert \int_{0}^{t} e^{(t-s)\Delta}  \P( \vN \cdot \F,  \vN \cdot \G)(s,\cdot) ds  \right\Vert_{T} \\   \nonumber 
&\leq & c (1+\sqrt{T})\left( \Vert \P( \vN \cdot \F) 
\Vert_{L^{2}_{t} \dot{H}^{-1}_{x}}  + \Vert \P( \vN \cdot \G)
\Vert_{L^{2}_{t} \dot{H}^{-1}_{x}}\right) \\
& \leq & c   (1+\sqrt{T})( \Vert \F \Vert_{L^{2}_{t} L^{2}_{x}}+ \Vert \G \Vert_{L^{2}_{t} L^{2}_{x}}).
\end{eqnarray}
Now, to study the bilinear terms recall the following estimate given in \cite{LR16} (Theorem $12.2$, page $352$):
\begin{Lemma} Let $ \vu, \vb \in E_T$. We have
\begin{equation*}
\begin{split}
    \ds{\left\Vert \int_{0}^{t} e^{(t-s)\Delta} \P((( \vu \ast \theta_{\varepsilon}) \cdot \nabla) \vb)(s,\cdot)ds \right\Vert_{T}  \leq c \sqrt{T} \varepsilon^{-3/2} \Vert \vu \Vert_{T} \Vert \vb \Vert_{T}}.    
\end{split}
\end{equation*}
\end{Lemma}	
Applying this lemma to each bilinear term in the equation (\ref{eq01}) we get 
\begin{equation}
\label{eq04}
    B((\vu,\vb), (\vv,\vc)) \leq c \sqrt{T}\varepsilon^{-3/2}  \Vert (\vu, \vb) \Vert_{T}\,  \Vert (\vv, \vc) \Vert_{T}.
\end{equation}
 
Once we have inequalities \eqref{eq02}, \eqref{eq03} and \eqref{eq04},  for a  time $0< T_0 <T_1 $ such that
\begin{equation*}
    T_0 = \min \left( T_1, \frac{c \varepsilon^3 }{(\| (\vu_0, \vb_0)\|_{L^2} + \|F\|_{L^2((0,T_1), L^2)}^2}) \right) ,
\end{equation*}
by the Picard's contraction principle, we obtain  $(\vu_\varepsilon, \vb_\varepsilon, \pe, \qe)$ a local solution of $(MHD_\varepsilon )$, where $ \vue, \vbe \in E_T $ and $\pe, \qe \in L^2((0,T), \dot H ^{-1}) +L^2((0,T), L^2)$. We can verify that this solution is unique.

To prove that $\pe \in L^4((0,T), L^{6/5}) +L^2((0,T), L^2)$, recall that 
$$\pe = \sum_{1 \leq i,j\leq 3} \mathcal{R}_i \mathcal{R}_j ((u_{\varepsilon, i}*\te) u_{\varepsilon,j} - (b_{\varepsilon,i}*\te) b _{\varepsilon,j} - F_{i,j} ),$$

As $\vue, \vbe \in E_T=L^{\infty}_{t}L^{2}_{x}\cap L^{2}_{t}\dot{H}^{1}_{x}$ then we have $\vue* \theta_\varepsilon, \vbe*\theta_\varepsilon \in E_T$ and thus we get $\vue,\vue* \theta_\varepsilon, \vbe, \vbe*\theta_\varepsilon \in L^{\infty}_{t}L^{2}_{x}\cap L^{2}_{t}L^{6}_{x}$. By interpolation we get $\vue*\theta_\varepsilon,\vbe*\theta_\varepsilon  \in L^{4}_{t}L^{3}_{x}$ and moreover, as $(\vue, \vbe)\in L^{\infty}_{t}L^{2}_{x}$ then by the H\"older inequalities,
$(\vue* \theta_\varepsilon)\otimes \vue, (\vbe* \theta_\varepsilon)\otimes \vbe \in L^{4}_{t}L^{6/5}_{x}$.
Thus, by the continuity of the Riesz transforms $\mathcal{R}_i$ on the Lebesgue spaces $L^p(\Rt)$ for $1<p<+\infty$ we have  $\sum_{1 \leq i,j\leq 3} \mathcal{R}_i \mathcal{R}_j ((u_{\varepsilon, i}*\te) u_{\varepsilon,j} - (b_{\varepsilon,i}*\te) b _{\varepsilon,j})\in L^{4}((0,T),L^{6/5})$.
Similarly we treat $\qe$.

Now, we prove that  $ (\vue, \vbe, \pe, \qe)$ is a global solution. We define the maximal existence time of the solution $\vu$ by
\begin{equation*}
T_{MAX} = \sup \{ 0< T \leq T_\infty \, : \,  \vu \in E_T  \}    
\end{equation*}
If $  T_{MAX} < T_{\infty}$ we take $ 0 < T< T_{MAX} < T_1 < T_{\infty}$, then  $(\vu, \vb)$ is a solution of ($GMHD_\varepsilon$) on $[0,T]$ and $(\vu, \vb)$ is a solution on $[T, T+ \delta]$, where
\begin{equation*}
    \delta = \min \left( T_1 - T , \frac{c \varepsilon^3 }{(\| (\vu(T), \vb(T))\|_{L^2} + \|F\|_{L^2((T,T_1), L^2)} )^2} \right),
\end{equation*}
which implies that $\lim_{T \to T_{MAX}^{-} } \|(\vue(T), \vbe(T) ) \|_{L^2} = + \infty$, however, we will see that it is not possible.

As $ ((\vb_{\varepsilon} * \te) \cdot \vN ) \vb_{\varepsilon} ) \vu_{_\varepsilon} = \vN \cdot ( \vb_{\varepsilon} \otimes ( \vb_{\varepsilon} * \te )) \vu_{_\varepsilon}$ belongs to $L^2((0,T), \dot H ^{-1})$, and the same for the other non linear terms, we can write
 \begin{eqnarray*}
 \frac{d}{dt} \Vert \vu_\varepsilon(t) \Vert^{2}_{L^2} &=& 2 \langle \partial_t \vu_\varepsilon (t), \vu_\varepsilon(t) \rangle_{\dot{H}^{-1}\times \dot{H}^{1}} \\
 &=& -2 \Vert \vN \vu_\varepsilon(t) \Vert^{2}_{L^2} +  2 \sum_{1\leq i,j \leq 3} \int   b_{\varepsilon,i} ( b_{\varepsilon,j} * \te ) \partial_i u_{_\varepsilon,j}  dx \\
 & & + 2 \sum_{1\leq i,j \leq 3} \int   F_{i,j}\partial_iu_{_\varepsilon,j}\    \, dx ,
 \end{eqnarray*}
  and 
 \begin{eqnarray*}
	\frac{d}{dt} \Vert \vb_\varepsilon(t) \Vert^{2}_{L^2} &=& 2 \langle \partial_t \vb_\varepsilon (t), \vb_\varepsilon(t) \rangle_{\dot{H}^{-1}\times \dot{H}^{1}} \\
	&=& -2 \Vert \vb_\varepsilon(t) \Vert^{2}_{\dot{H}^{1}} +  2 \sum_{1\leq i,j \leq 3} \int   u_{\varepsilon,i} ( b_{\varepsilon,j} * \te ) \partial_i b_{_\varepsilon,j}  dx  \\
	& & + 2 \sum_{1\leq i,j \leq 3} \int   G_{i,j}\partial_iu_{_\varepsilon,j}\    \, dx.
\end{eqnarray*}  
where we have used the fact that 
\begin{align*}
\label{d}
    \int ( (\vu_\varepsilon* \theta ) \cdot \nabla) \vb_\varepsilon  \cdot \vb_\varepsilon \, dx &= \int \sum_{1\leq i,j \leq 3} ( (u_{j, _\varepsilon} *\theta ) \partial_j b_{i,_\varepsilon}) b_{ i, _\varepsilon}\, dx \nonumber \\
    &= -\frac{1}{2} \int  (\vu_{ \varepsilon} *\theta ) \cdot \vN ( |\vb_\varepsilon|^2 ) \, dx \\ &=-\frac{1}{2} \int  \nabla\cdot(\vue *\theta_\varepsilon)  |\vb_\varepsilon|^2  dx 
    = 0.
\end{align*}

Then, an integration by parts gives
$$ \sum_{1\leq i,j \leq 3} \int   u_{\varepsilon,i} ( b_{\varepsilon,j} * \te ) \partial_i b_{\varepsilon,j}  dx  = -  \sum_{1\leq i,j \leq 3} \int   b_{\varepsilon,i} ( b_{\varepsilon,j} * \te ) \partial_i u_{\varepsilon,j}  dx,$$ 
so we have
\begin{eqnarray*}
 \frac{d}{dt} ( \Vert \vu_\varepsilon(t) \Vert^{2}_{L^2}+ \Vert \vb_\varepsilon(t) \Vert^{2}_{L^2})&=&-2 (\Vert \vN \vu_\varepsilon(t) \Vert^{2}_{L^2} + \Vert \vN \vb_\varepsilon(t) \Vert^{2}_{L^2} )	\\
 & & + 2 \sum_{1 \leq i,j \leq 3} (   \int   F_{i,j}\partial_iu_{j}\    \, dx\, ds + \int   G_{i,j}\partial_ib_{j}\    \, dx\, ds ) . \\ 
\end{eqnarray*}	
By integrating on the time interval $[0,T]$ we obtain the control (\ref{eneq1})  which implies by Gr\"onwall inequality that $\|(\vue, \vbe)(T) ) \|_{L^2}$ does not converges to $+\infty$ when $T$ go to $T_{MAX}$ if $T_{MAX}< T_\infty$, hence the solution is defined on $[0, T_\infty)$. Finally,  remark that we can write
\begin{align*}
    \vN \cdot ((\vbe \cdot \vue) (\vbe*\theta_\varepsilon)) &= \vN (\vbe \cdot \vue) \cdot (\vbe*\theta_\varepsilon) \\
    &=((\vbe*\theta_\varepsilon) \cdot \vN) \vbe \cdot \vue + ((\vbe*\theta_\varepsilon) \cdot \vN) \vue )  \cdot \vbe 
\end{align*}
so that
\begin{align*}
    \partial_t(\frac {\vert\vu_\varepsilon \vert^2}2)=&\Delta(\frac {\vert\vu_\varepsilon \vert^2}2)-\vert\vN \vu_\varepsilon \vert^2- \vN\cdot\left( \frac{\vert\vu_\varepsilon \vert^2}2 (\vu_\varepsilon * \theta_\varepsilon ) + \pe \vue \right) \\
    &+ \vN \cdot ((\vu_\varepsilon \cdot \vb_\varepsilon) (\vb_\varepsilon * \theta_\varepsilon ))  - ((\vbe*\theta_\varepsilon) \cdot \vN) \cdot \vue ) \vbe  + \vu_\varepsilon \cdot(\vN\cdot\mathbb{F}),
\end{align*}
similarly we find
\begin{align*}
    \partial_t(\frac {\vert\vb_\varepsilon \vert^2}2)=&\Delta(\frac {\vert\vb_\varepsilon \vert^2  }2)-\vert\vN \vb_\varepsilon \vert^2 -\vN\cdot\left(\frac{\vert\vb_\varepsilon \vert^2}2(\vu_\varepsilon * \theta_\varepsilon ) + \qe \vbe \right) \\
    &+ ((\vbe*\theta_\varepsilon) \cdot \vN) \cdot \vue ) \vbe + \vbe  \cdot(\vN\cdot\mathbb{G}).
\end{align*}
By adding these equations we obtain the energy equality \eqref{eneqpo}.
\Endproof \\ 

We can observe that our approximated system need to consider an non-zero term $q_\varepsilon$ even if $G=0$. As we have seen it is not the case when we let $\epsilon$ tends to $0$ and then we obtain the (MHDG) system.

\end{document}